\documentclass[10pt]{article}
\usepackage{tocloft}
\usepackage{amsmath}
\usepackage{amsfonts}
\usepackage{amssymb}
\usepackage{amsthm}
\usepackage{bbm}
\renewcommand{\sqrt}[1]{\left( #1 \right)^\frac12}
\newcommand{\concat}{\frown}

\newcommand{\e}[0]{\mathrm{e}}
\newcommand{\E}[0]{\mathbb{E}}
\newcommand{\T}[0]{\mathbb{T}}
\newcommand{\Ber}[0]{\mathrm{Ber}}
\newcommand{\pv}{\mathrm{P.V. }}
\renewcommand{\i}[0]{i}
\renewcommand{\d}[0]{\mathrm{d}}
\newcommand{\sgn}[0]{\mathrm{sgn}}

\newcommand{\indnorm}[1]{\left\|#1\right\|_{\mathrm{ind}}}
\newcommand{\Lnorm}[1]{\left\|#1\right\|_{L^1\left(\ell^2\right)}}
\usepackage[utf8]{inputenc}
\author{Maciej Rzeszut, Michał Wojciechowski}

\title{Independent sums of $H^1_n(\T)$ and $H^1_n(\delta)$}

\newtheorem{thm}{Theorem}[section]
\newtheorem{lem}[thm]{Lemma}
\newtheorem{deff}[thm]{Definition}
\newtheorem{cor}[thm]{Corollary}

\begin{document}
\maketitle

\begin{abstract}
We construct a new idempotent Fourier multiplier on the Hardy space on the bidisc, which could not be obtained by applying known one dimentional results. The main tool is a new $L^1$ equivalent of the Stein martingale inequality which holds for a special filtration of periodic subsets of $\mathbb{T}$ with some restrictions on the functions involved. We also identify the isomorphic type of the range of the associated operator as the independent sum of dyadic $H^1_n$, which is known to be a complemented and invariant subspace of dyadic $H^1$. 
\end{abstract}

\tableofcontents
\newpage
\section*{Introduction}
\addcontentsline{toc}{section}{Introduction}
The basic notion of this paper is one of the possible ways to define a norm on a direct sum of a sequence of certain Banach spaces. 

\begin{deff} Let $X_n$ be a sequence of subspaces of $L^1\left(\Omega,\mathcal{F},\mu,B\right)$, where $\mu$ is a probability measure. By the independent sum of $X_1,\ldots,X_n,\ldots$ we mean the space $\left(\bigoplus_{n=1}^\infty X_n\right)_{\mathrm{ind}}$ of sequences $\left(f_1,\ldots,f_n,\ldots\right)$ such that $f_n\in X_n$ and the norm \[\left\|\left(f_1,\ldots,f_n,\ldots\right)\right\|_{\mathrm{ind}}= \int_{\Omega^\infty}\left(\sum_{n=1}^\infty \left\|f_n\left(\omega_n\right)\right\|_B^2\right)^\frac12 \d \left(\mu\left(\omega_1\right)\otimes\ldots\otimes\mu\left(\omega_n\right)\otimes\ldots\right)\] is finite, equipped with this norm.\end{deff}

Thus, the independent sum of subspaces of $L^1\left(\Omega,B\right)$ can be identified with a subspace of $L^1\left(\Omega^\infty,\ell^2(B)\right)$. In particular, it is a Banach space. Care has to be taken, since contrary to $\ell^p$- and $c_0$-direct sums, uniformly bounded operators acting on $X_n$ do not necessarily induce a bounded operator on $\left(\bigoplus X_n\right)_{\mathrm{ind}}$ and thus there is no apparent reason for the isomorphism class of the independent sum to depend only on isomorphism classes of the summands. Remedying this issue will be one of the main concerns of Section 5.  \par 
The paper is organised as follows. In Section 1 we develop tools needed in the following sections to handle certain special cases of the Stein martingale inequality for $p=1$, in particular an improved version of the so-called non-linear telescoping lemma, originally due to Bourgain \cite{bourgtelesc} and later put to use by M\"{u}ller \cite{mullerdecomp}. In Section 2 we provide an example of a not yet known complemented and invariant subspace of the Hardy space on the bidisc and prove that it is isomorphic to the independent sum $\left(\bigoplus H^1_n(\T)\right)_{\mathrm{ind}}$. In Section 3 we introduce the space $\left(\sum H^1_n(\delta)\right)_{\mathrm{ind}}\subset H^1(\delta)$ (which is isometric to the independent sum $\left(\bigoplus \iota\left(H^1_n(\delta)\right)\right)_{\mathrm{ind}}$, where $\iota:H^1(\delta)\to L^1\left(\ell^2\right)$ is the canonical isometric embedding) and give a new proof to the known result of M\"{u}ller and Schechtman \cite{smallmuller}, stating that this space is complemented in $H^1(\delta)$. In Section 4 we prove that the orthogonal projection $P:L^2\left(\ell^2\right)\to \iota\left(H^2(\delta)\right)$ and the operator $\iota^{-1}P:L^2\left(\ell^2\right)\to H^2(\delta)$ are Calder\'{o}n-Zygmund operators and extend via principal value to weak type operators on respctive $L^1$ spaces, thus by a theorem of Bourgain \cite {bourginterp} yielding the $K$-closedness of the couple $\left(\iota\left(H^1(\delta)\right),\iota\left(H^2(\delta)\right)\right)$ in $\left(L^1\left(\ell^2\right),L^2\left(\ell^2\right)\right)$. In Section 5, using a theorem of Johnson and Schechtman \cite{JSch}, we prove a general theorem about isomorphisms of independent sums and provide two ways to deduce the isomorphism announced in title from it: one based on Wojtaszczyk's \cite{smallwojt} construction of an isomorphism between $H^1(\T)$ and $H^1(\delta)$ and the $K$-closedness, and the other based on later Meyer's construction \cite{meyer}, without the use of $K$-closedness. In Section 6 we conclude and point out some open questions raised by the proved results.\par 
We would like to thank P. F. X. M\"{u}ller for many helpful comments and suggestions.
\section{Non-linear telescoping and related results}

In the whole paper, we will denote increasing filtrations by $\left(\mathcal{F}_n\right)$ and decreasing filtrations by $\left(\mathcal{F}^*_n\right)$. To shorten the notiation we will write $\E_n f$ instead of $\E\left(f\mid \mathcal{F}_n\right)$ and $\E^*_n f$ instead of $\E\left(f\mid \mathcal{F}^*_n\right)$. \par 
By the signs $\gtrsim$ $, \lesssim$, $\sim$ we mean respectively $\geq,\leq,=$ up to a constant. Let us recall the classical Stein martingale inequality \cite{Stein}. 
\begin{thm}\label{classicstein}Let $\left(\mathcal{F}^*_k\right)_{k=1}^n$ be a decreasing filtration on a probability space $\left(\Omega,\mathcal{F},\mu\right)$ and $f_1,\ldots,f_n$ be integrable $\mathcal{F}$-measurable functions. Then for $1<p<\infty$,
\[\E\left(\sum_{k=1}^n\left|f_k\right|^2\right)^\frac{p}{2}\gtrsim \E\left(\sum_{k=1}^n\left|\E^*_k f_k\right|^2\right)^\frac{p}{2}\] with a constant depending only on $p$. \end{thm}
Obviously, the order of the filtration has no influence on the inequality, but for reasons of consistency with subsequent considerations we choose to use the decreasing order. This inequality is false for $p=1$ in general, but under certain additional assumptions it remains true for $p=1$. One of results of this type is the following inequality of Lepingle \cite{Lep}.
\begin{thm} Let $\left(\mathcal{F}^*_k\right)_{k=1}^n$ be a decreasing filtration on a probability space $\left(\Omega,\mathcal{F},\mu\right)$ and $f_1,\ldots,f_n$ be integrable functions such that $f_k$ is $\mathcal{F}^*_{k-1}$-measurable. Then 
\[\E\left(\sum_{k=1}^n\left|f_k\right|^2\right)^\frac{1}{2}\geq \frac{1}{2}\E\left(\sum_{k=1}^n\left|\E^*_k f_k\right|^2\right)^\frac{1}{2}.\]\end{thm}
The following lemma will be used to estimate the left-hand side in Stein-type inequalities from below. 
\begin{lem}\label{inf}Let $\left(\mathcal{F}^*_k\right)_{k=1}^n$ be a decreasing filtration and $\left(\varphi_k\right)_{k=1}^n$ be an adapted sequence of integrable functions. Define $\lambda_0=0$, $\lambda_k=\E^*_{k+1}\left(\left|\varphi_k\right|^2+\lambda_{k-1}^2\right)^\frac12$ for $k=1,\ldots,n$. Then 
\[\inf_{\E^*_k f_k=\varphi_k, k=1,\ldots,n} \E\left(\sum_{k=1}^n\left|f_k\right|^2\right)^\frac{1}{2}=\E\lambda_n.\]\end{lem}
\textit{Proof.} We will prove the '$\geq$' inequality first. Note that for $k\geq 1$, 
\[\E\left(\lambda _{k-1}^2+\sum_{j=k}^n \left|\E_k^* f_j\right|^2\right)^\frac{1}{2}= \E\left( \lambda_{k-1}^2+\left|\E^*_kf_k\right|^2 + \sum_{j=k+1}^n \left|\E^*_k f_j\right|^2\right)^\frac{1}{2}\geq\]
\[\geq  \E\left( \left(\E^*_{k+1}\left(\lambda_{k-1}^2+\left|\E^*_kf_k\right|^2\right)^\frac12\right)^2 + \sum_{j=k+1}^n \left|\E^*_{k+1}\E^*_k f_j\right|^2\right)^\frac12= \]\[=\E\left(\lambda_k^2+\sum_{j=k+1}^n\left|\E^*_{k+1}f_j\right|^2\right)^\frac{1}{2}\]
and thus by induction, $\E\left(\sum_{j=1}^n \left|f_j\right|^2\right)^{\frac{1}{2}}\geq \E\left(\lambda_k^2+\sum_{j=k+1}^n\left|\E^*_{k+1}f_j\right|^2\right)^\frac{1}{2}$. In particular, $\E\left(\sum_{j=1}^n \left|f_j\right|^2\right)^{\frac{1}{2}}\geq \E\lambda_n$. \par
Now we will prove the '$\leq$' inequality. Assume first that $\varphi_k$ satisfy $\frac{1}{M}<\left|\varphi_k\right|<M$ for some constant $M$ and denote $\Lambda_k=\left(\left|\varphi_k\right|^2+\left(\E^*_k\Lambda_{k-1}\right)^2\right)^\frac12$, $\Lambda_0=0$. Obviously $\lambda_k=\E^*_{k+1}\Lambda_k$. Take \[f_k=\varphi_k\prod_{j=1}^{k-1}\frac{\Lambda_j}{\lambda_j}.\] 
This choice is legal, since the conditions imposed on $\varphi_k$ ensure that $f_k$ is integrable and for any $j=2,\ldots,k$ we have \[\E^*_{k}\left(\varphi_k\prod_{i=j-1}^{k-1}\frac{\Lambda_{i}}{\lambda_i}\right)=\E^*_k\E^*_j\left(\frac{\Lambda_{j-1}}{\lambda_{j-1}}\varphi_k\prod_{i=j}^{k-1}\frac{\Lambda_i}{\lambda_i}\right)=\]\[=\E^*_k\left(\left(\E^*_j\frac{\Lambda_{i-1}}{\lambda_{i-1}}\right)\varphi_k\prod_{i=j}^{k-1}\frac{\Lambda_i}{\lambda_i}\right)= \E^*_k\left(\varphi_k\prod_{i=j}^{k-1}\frac{\Lambda_i}{\lambda_i}\right)\] (by the fact that $\E\left(uv\mid \mathcal{F}\right)=v\E\left(u\mid\mathcal{F}\right)$ for $v$ bounded and $\mathcal{F}$-measurable) and thus by induction $\E^*_k f_k=\E^*_{k}\left(\varphi_k\prod_{i=j-1}^{k-1}\frac{\Lambda_{i}}{\lambda_i}\right)$, in particular $\E^*_k f_k=\varphi_k$. For $k\geq 2$, \[\E\left(\Lambda_{k-1}^2+\sum_{j=k}^n\left|\varphi_j\prod_{i=k-1}^{j-1}\frac{\Lambda_i}{\lambda_i}\right|^2\right)^\frac{1}{2}= \E\left(\frac{\Lambda_{k-1}}{\lambda_{k-1}}\left(\lambda_{k-1}^2+ \sum_{j=k}^n\left|\varphi_j\prod_{i=k}^{j-1}\frac{\Lambda_i}{\lambda_i}\right|^2\right)^\frac{1}{2}\right)=\] \[\E\left(\lambda_{k-1}^2+ \sum_{j=k}^n\left|\varphi_j\prod_{i=k}^{j-1}\frac{\Lambda_i}{\lambda_i}\right|^2\right)^\frac{1}{2}=\E\left(\lambda_{k-1}^2+\left|\varphi_{k}\right|^2+ \sum_{j=k+1}^n\left|\varphi_j\prod_{i=k}^{j-1}\frac{\Lambda_i}{\lambda_i}\right|^2\right)^\frac{1}{2}=\] \[\E\left(\Lambda_k^2+\sum_{j=k+1}^n\left|\varphi_j\prod_{i=k}^{j-1}\frac{\Lambda_i}{\lambda_i}\right|^2\right)^\frac{1}{2}.\] Since $\Lambda_1=\left|\varphi_1\right|$, we have by induction \[\E\left(\sum_{j=1}^n\left|f_j\right|^2\right)^\frac{1}{2}= \E\left(\Lambda_{k-1}^2+\sum_{j=k}^n\left|\varphi_j\prod_{i=k-1}^{j-1}\frac{\Lambda_i}{\lambda_i}\right|^2\right)^\frac{1}{2}\]. In particular, $\E\left(\sum_{j=1}^n\left|f_j\right|^2\right)^\frac{1}{2}=\E\Lambda_n=\E\lambda_n$. Now let $\varphi_k$ be arbitrary and define \[\varphi_k^{(M)}(\omega)=\left\{\begin{array}{lll}M\sgn\left(\varphi_k(\omega)\right)&\ \mathrm{if}\ &\left|\varphi_k(\omega)\right|>M\\ \frac{1}{M}\sgn\left(\varphi_k(\omega)\right)&\ \mathrm{if}\ &\left|\varphi_k(\omega)\right|<\frac{1}{M}\\ \varphi_k(\omega)&\ \mathrm{otherwise.}\ &\ \end{array}\right.\] Obviously $\varphi_k^{(M)}$ are still $\mathcal{F}_k^*$-measurable. We can construct $\lambda^{(M)}_k$, $\Lambda^{(M)}_k$ and $f^{(M)}_k$ out of $\varphi^{(M)}_k$ just as we did with $\varphi_k$. Now take $f_k=f_k^{(M)}+\varphi_k-\varphi_k^{(M)}$. We have $\E^*_k f_k= \varphi_k^{(M)}+\varphi_k-\varphi_k^{(M)}=\varphi_k$ and by the previous calculation \[\E\left(\sum_{k=1}^n\left|f_k\right|^2\right)^\frac{1}{2}\leq \E\left(\sum_{k=1}^n\left|f_k^{(M)}\right|^2\right)^\frac{1}{2}+ \E\left(\sum_{k=1}^n\left|f_k-f_k^{(M)}\right|^2\right)^\frac{1}{2}\leq\]\[ \E\lambda_n^{(M)}+\sum_{k=1}^n\E\left|f_k-f^{(M)}_k\right|\] which tends to $\E\lambda_n$ when $M\to \infty$ by the conditional dominated convergence theorem, thus proving the desired inequality. \par 
We will be mainly concerned with the '$\geq$' inequality in Lemma \ref{inf}. The following theorem, which up to rephrasing appears as an unproved remark in \cite{BRSch}, shows an application of the '$\leq$' inequality. It has been communicated by P. F. X. M\"{u}ller that it can be proved by an explicit construction, although we find the subsequent proof more fun.

\begin{thm}\label{drzewa} Let $f_0,\ldots,f_n\in L^1\left(\{0,1\}^n\right)$ and $\mathcal{F}^*_k$ be the product of the trivial sigma-algebra of subsets of $\{0,1\}^k$ and the full sigma-algebra of subsets of $\{0,1\}^{n-k}$ for $k=0,\ldots,n$. Then there is no universal constant for the inequality
\[\E\left(\sum_{k=0}^n \left|f_k\right|^2\right)^{\frac{1}{2}}\gtrsim \E\left(\sum_{k=0}^n \left|\E^*_k f_k\right|^2\right)^{\frac{1}{2}}\] under the assumption $\E^*_{k+1}f_k=0$ to hold (the choice of $\mathcal{F}^*_{n+1}$ is irrelevant - we may take for instance $\mathcal{F}^*_{n+1}=\mathcal{F}^*_{n}$).\end{thm}
\textit{Proof.} Measurability with respect to $\mathcal{F}^*_k$ is equivalent to being dependent only on last $n-k$ variables. More precisely,
\[\left(\E^*_k f\right)\left(x_1,\ldots,x_n\right)=\frac{1}{2^k}\sum_{y_1,\ldots,y_k\in \left\{0,1\right\}} f\left(y_1,\ldots,y_k,x_{k+1},\ldots,x_n\right).\]
Let $\varphi_k=\E^*_k f_k$ and $\lambda_k$ be defined by $\lambda_{-1}=0$, $\lambda_{k}=\E^*_{k+1}\left(\left|\varphi_k\right|^2+\lambda_{k-1}^2\right)^\frac12$ for $k=0,\ldots,n$. The condition $\E^*_{k+1}f_k=0$ is equivalent to $\E^*_{k+1} \varphi_k=0$. We may now apply Lemma \ref{inf} and reduce our problem to proving the the non-inequality $\E\lambda_n\not\gtrsim\E\left(\sum\left|\varphi_k\right|^2\right)^\frac{1}{2}$, in which only $\varphi_k$ appear.\par 

Let $T_n$ (respectively $T_n^*$) be the set of scalar-valued functions on the full binary tree of height $n$ (respectively the full binary tree of height $n$ with an additional vertex attached to the root) and denote $T=\bigcup T_n$, $T^*=\bigcup T_n^*$. For $U_1,U_2\in T_n$, $a\in \mathbb{R}$ we can naturally define the concatenations $a\concat U_1\in T_n^*$, $U_1\concat a \concat U_2\in T_{n+1}$. Let us define the norm $|||\cdot|||$ on $T$ recursively. For $U\in T_0$, $|||U|||=|u|$, where $u$ is the value of $U$ in the single element of the domain and for $U\in T$, $U=U_1\concat a \concat U_2$ we put \[||| U|||=\frac{1}{2}\left(\left(a^2+|||U_1|||^2\right)^\frac12+\left(a^2+|||U_2|||^2\right)^\frac12\right).\]
We can extend this definition to $U\in T^*$ by $|||U|||=\left(a^2+|||U_1|||^2\right)^\frac12$ for $U=a\concat U_1$. An alternate definition is as follows. Let $T_n\cup T_n^*\ni U\mapsto U'\in T_{n-1}\cup T_{n-1}^*$ be an operation that removes the two downmost rows and puts $\frac{1}{2}\left(\left(a^2+b^2\right)^\frac12+\left(a^2+c^2\right)^\frac12\right)$ in the place of any subtree $b\concat a\concat c$ consisting of two leaves $b,c$ and their father $a$. Then $|||U|||=U^{\prime\prime\cdots\prime}$. Moreover, we define the norm $||\cdot||$ on $T\cup T^*$ by
\[\|U\|=\frac{1}{2^n} \sum_{x}\left(U\left(x\right)^2+ U\left(\tilde{x}\right)^2+U\left(\tilde{\tilde{x}}\right)^2+\ldots\right)^\frac12,\]
where $n$ is the height of $U$, summation is taken over the leaves and $\tilde{x}$ denotes the father of $x$. \par 
Let us return to the previous picture. Since all conditions imposed on $\varphi_k$ are $\mathcal{F}^*_k$-measurability (equivalent to dependence only on $x_{k+1},\ldots,x_{n}$) and $\E^*_{k+1}\varphi_k=0$, the dependence on $x_{k+1}$ is contained in the factor $(-1)^{x_{k+1}}$ and thus $\left|\varphi_k\left(x_1,\ldots,x_n\right)\right|=g_k\left(x_{k+2},\ldots,x_n\right)$ for some nonnegative $g_k$. We can think of the sequence $\left(g_0,\ldots,g_n\right)$ as an element $G$ of $T^*_{n-1}$. Obviously $\|G\|= \E\left(\sum\left|\varphi_{k}\right|^2\right)^{\frac12}$. By $\lambda_{k}=\E^*_{k+1}\left(g_k^2+\lambda_{k-1}^2\right)^\frac12$ and induction it follows that the rows of $G^{\prime(k)}$ are $\lambda_k,g_{k+1},\ldots,g_n$ and thus $|||G|||=\lambda_n=\E\lambda_n$ (since $\lambda_n$ is constant). The only thing left to prove is that there is no constant independent on $n$ for the inequality $|||G|||\gtrsim ||G||$ to hold for $G\in T^*$.\par
Let \[f(t)=\inf_{\substack{G\in T\\\|1\concat G\|=t}}|||G|||\] for $t\geq 1$. Limiting the range of infimum to $G=U\concat 0\concat V$, for any $x,y\geq 1$ such that $\frac{x+y}{2}=t$ we have 
\[\aligned f(t)&\leq \inf_{\substack{U,V\in T\\ \|1\concat U\|=x\\ \|1\concat V\|=y}}|||U\concat 0\concat V|||\\&= \frac{1}{2}\left(f(x)+f(y)\right).\endaligned\] It is routine to check that $f$ is measurable, so by a theorem of Sierpi\'{n}ski \cite{Sierp} $f$ is convex and thus it is a.e. differentiable and equal to the integral of its derivative. Now let $1\leq a\leq t$. Then 
\[\aligned f(t)&\leq \inf_{\substack{G\in T\\ \left\| 1\concat\left(G\concat \sqrt{a^2-1}\concat G\right) \right\|=t}}|||G\concat \sqrt{a^2-1}\concat G|||\\ &=\inf_{\|a\concat G\|=t}\sqrt{a^2-1+|||G|||^2}\\&=\inf_{\|1\concat G\|=\frac{t}{a}} \sqrt{a^2-1+a^2|||G|||^2}\\&= \sqrt{a^2-1+a^2f\left(\frac{t}{a}\right)^2}. \endaligned\] Thus the function $t\mapsto \frac{1+f(t)^2}{t^2}$ is nonincreasing. For computational convenience let us change the variable from $t$ to $t-1$. Then $f(0)=0$ and $\frac{1+f(t)^2}{(1+t)^2}$ is nonincreasing. The latter gives \[(t+1)f(t)f'(t)\leq 1+f(t)^2.\] By convexity, $f(t)=\int_{0}^{t}f'(\tau)\d \tau\leq tf'(t)$ a.e. Combining these two inequalites we get $f(t)^2\leq t$, which together with convexity gives $f(t)=0$ for almost all $t$. Taking arbitrarily large $t$ and $G$ such that $\|1\concat G\|=t$, $|||G|||<\varepsilon$, we can make $\frac{|||1\concat G|||}{\|1\concat G\|}\leq \frac{1+\varepsilon}{t}$ arbitrarily small, which ends the proof.\par 

We may now proceed to the announced non-linear telescoping lemma. We will prove a variation of a lemma extracted by M\"{u}ller in \cite{mullerdecomp} from the work of Bourgain \cite{bourgtelesc}. Although we get a (much) worse estimate in the general case, we get the estimates for the $\|\cdot\|_{\mathrm{ind}}$ norm. 

\begin{lem}\label{telescoping}Let $\lambda_k$ and $\varphi_k$ for $k=1,\ldots,n$ be nonnnegative integrable random variables. Put $\lambda_0=0$. Then\\
(i) if $\lambda_k=\E\sqrt{\varphi_k^2+\lambda_{k-1}^2}$, then \[\max\left(\left\|\left(\varphi_{k}\right)_{k=1}^n\right\|_{L^1\left(\ell^2\right)}, \left\|\left(\varphi_{k}\right)_{k=1}^n\right\|_{\mathrm{ind}}\right)\leq 2\lambda_n\]
(ii) if $\E\lambda_k\geq \E\sqrt{\varphi_k^2+\lambda_{k-1}^2}$ and $\lambda_k$ are bounded, then \[\max\left(\left\|\left(\varphi_{k}\right)_{k=1}^n\right\|_{L^1\left(\ell^2\right)}, \left\|\left(\varphi_{k}\right)_{k=1}^n\right\|_{\mathrm{ind}}\right)\leq \left(1+2^\frac12\right)\sqrt{\E\lambda_n}\sqrt{\max_{k=1,\ldots,n}\sup \lambda_k}.\]\end{lem}
\textit{Proof.} Denote $g_k=\sqrt{\varphi_k^2+\lambda_{k-1}^2}-\lambda_{k-1}$. Then $\varphi_k^2=\left(g_k+\lambda_{k-1}\right)^2-\lambda_{k-1}^2=g_k^2+2g_k\lambda_{k-1}$ and thus $\varphi_k\leq g_k+\sqrt{2g_k\lambda_{k-1}}$. Note that both of the norms $\left\|\cdot\right\|_{L^1\left(\ell^2\right)}$ and $\left\|\cdot\right\|_{\mathrm{ind}}$ are not greater than any of the norms $\left\|\cdot\right\|_{L^1\left(\ell^1\right)}$ and $\left\|\cdot\right\|_{L^2\left(\ell^2\right)}$. Thus 
\[\max\left(\left\|\left(\varphi_{k}\right)_{k=1}^n\right\|_{L^1\left(\ell^2\right)}, \left\|\left(\varphi_{k}\right)_{k=1}^n\right\|_{\mathrm{ind}}\right)\leq \left\|\left(g_k\right)_{k=1}^n\right\|_{L^1\left(\ell^1\right)}+ \left\|\left(\sqrt{2g_k\lambda_{k-1}}\right)_{k=1}^n\right\|_{L^2\left(\ell^2\right)}\]
\[\leq\sum_{k=1}^n\E g_k+\sqrt{\sum_{k=1}^n\E \left(2g_k\lambda_{k-1}\right)}.\]
In order to prove (i), we notice that $\E g_k=\lambda_k-\lambda_{k-1}$ and exploit the telescoping nature of sums of $\E g_k$ to get \[\sum_{k=1}^n\E g_k+\sqrt{\sum_{k=1}^n\E \left(2g_k\lambda_{k-1}\right)}=\sum_{k=1}^n\left(\lambda_k-\lambda_{k-1}\right)+\sqrt{2\sum_{k=1}^n\lambda_{k-1}\E g_k}=\]\[\lambda_n+\sqrt{2\sum_{k=1}^n\E g_k\sum_{j=1}^{k-1}\E g_j}\leq\lambda_n+\sqrt{\sum_{j,k=1}^n\E g_j\E g_k}=2\lambda_n.\]
Now we will prove (ii). In this case $\E g_k\leq \E\lambda_k-\E\lambda_{k-1}$. 
\[\sum_{k=1}^n\E g_k+\sqrt{\sum_{k=1}^n\E \left(2g_k\lambda_{k-1}\right)}\leq \E\lambda_n+\sqrt{2\max_{k\leq n-1}\sup\lambda_k\sum_{k=1}^n\E g_k}\leq \]\[\E\lambda_n+2^\frac{1}{2}\sqrt{\max_{k\leq n-1}\sup \lambda_k}\sqrt{\sum_{k=1}^n \E g_k}\leq \left(1+2^\frac12\right)\sqrt{\E\lambda_n}\sqrt{\max_{k\leq n}\sup \lambda_k}.\]
The following observation shows that the estimates for the $\|\cdot\|_{L^1\left(\ell^2\right)}$ norms in Lemma \ref{telescoping} are, up to a constant, redundant.
\begin{cor}\label{2indgeqL1l2}Let $\varphi_1,\ldots,\varphi_n$ be integrable random variables and $\psi_1,\ldots,\psi_n$ be independent random variables such that for any $k$, the distributions of $\psi_k$ and $\varphi_k$ are the same. Then \[\E\left(\sum\left|\varphi_k\right|^2\right)^\frac{1}{2}\leq 2\E\sqrt{\sum\left|\psi_k\right|^2}.\]
In other words, $\left\|\left(\varphi_k\right)\right\|_{L^1\left(\ell^2\right)}\leq 2 \left\|\left(\varphi_k\right)\right\|_{\mathrm{ind}}$. \end{cor}
\textit{Proof.} Let $\lambda_k=\E\sqrt{\left|\varphi_k\right|^2+\lambda_{k-1}^2}$. Since $\lambda_k$ depends only on $\lambda_{k-1}$ and the distribution of $\varphi_k$, we have $\lambda_k=\E\sqrt{\left|\psi_k\right|^2+\lambda_{k-1}^2}$. Let $\mathcal{F}^*_k=\sigma\left(\psi_{k},\ldots,\psi_n\right)$. Then also $\lambda_k=\E^*_{k+1}\sqrt{\left|\psi_k\right|^2+\lambda_{k-1}}$, since $\psi_k$ and $\mathcal{F}^*_{k+1}$ are independent. Thus \[\E\left(\sum\left|\varphi_k\right|^2\right)^\frac{1}{2}\leq 2\lambda_n =2\E\lambda_n \leq 2\E\sqrt{\sum\left|\psi_k\right|^2},\] where the first inequality is an application of part (i) of Lemma \ref{telescoping} to $\left|\varphi_k\right|$ and the second is an application of Lemma \ref{inf} to $\psi_k$ and $\mathcal{F}^*_k$.\par 
The function $t\mapsto t^\frac12$ in Corollary \ref{2indgeqL1l2} can be replaced by an arbitrary concave and nonincreasing function with value $0$ in $0$ and the constant $2$ can be improved to $\frac{\e}{\e-1}$, but this fact is of no interest for this paper. We note one more consequence of the non-linear telescoping lemma. 
\begin{cor}\label{duallep}Let be a decreasing filtration. Let $f_1,\ldots,f_n$ be integrable functions such that $\sigma\left(f_k\right)\cap\mathcal{F}^*_{k+1}$ is trivial. Then \[\E\left(\sum_{k=1}^n\left|f_k\right|^2\right)^\frac{1}{2}\geq \frac{1}{2}\E\left(\sum_{k=1}^n\left|\E^*_k f_k\right|^2\right)^\frac{1}{2}.\]\end{cor}
\textit{Proof.} As usual, we take $\lambda_0=0$, $\lambda_{k}=\E^*_{k+1}\sqrt{\left|\E^*_kf_k\right|^2+\lambda_{k-1}^2}$ and note that it follows by induction that $\lambda_k$ is constant, since if $\lambda_{k-1}$ is constant, then $\sqrt{\left|\E^*_kf_k\right|^2+\lambda_{k-1}^2}$ is $\sigma\left(f_k\right)$-measurable and thus $\lambda_k$ is $\sigma\left(f_k\right)\cap\mathcal{F}^*_{k+1}$-measurable, since conditional expectations commute. Therefore \[\E\left(\sum_{k=1}^n\left|f_k\right|^2\right)^\frac{1}{2}\geq \E\lambda_n= \lambda_n\geq \frac{1}{2}\E\left(\sum_{k=1}^n\left|\E^*_k f_k\right|^2\right)^\frac{1}{2},\]
where the first inequality is an application of Lemma \ref{inf} to $f_k$ and $\mathcal{F}^*_k$ and the second is an application of part (i) of Lemma \ref{telescoping} to $\E^*_k f_k$. \par 

\section{Independent sum of $H^1_n(\T)$}

We first recall basic facts concerning invariant operators. Let $G$ be a compact abelian metric group and $\Gamma$ be its dual group. For $g\in G$ we define the shift $\tau_g:L^1(G)\to L^1(G)$ by $\tau_g(f)\left(g_0\right)=f\left(g_0-g\right)$. The functions $G\ni g\mapsto \tau_g(f)\in L^1(G)$ are known to be continuous for any $f$. Let $X$ be a shift-invariant subspace of $L^1(G)$, which in our case is equivalent to being spanned by a set of characters, say $\Gamma_0\subset \Gamma$. Any bounded invariant operator $T:X\to X$ corresponds to its Fourier multiplier, i.e the function $\left(t_\gamma\right)_{\gamma\in\Gamma_0}$ such that $\widehat{Tf}(\gamma)=t_\gamma \widehat{f}(\gamma)$ for $\gamma \in \Gamma_0$. Since for any regular measure $\mu$ of bounded variation, the multiplier of the operator $L^1(G)\ni f\mapsto \mu \ast f\in L^1(G)$ is the sequence $\left(\widehat{\mu}_\gamma\right)_{\gamma\in\Gamma}$, by an abuse of notation, we will write $\widehat{T}(\gamma)=t_\gamma$ in general case. \par
An invariant operator $T:X\to X$ is called idempotent if $T^2=T$ or, equivalently, $\widehat{T}\left(\gamma\right)\in\{0,1\}$. Bounded invariant operators are in a one-to-one correspondence with invariant and complemented subsapces of $X$. Indeed, if $T$ is invariant and idempotent, then the range of $T$ is $\overline{\mathrm{span}}\left(\mathrm{supp} \widehat{T}\right)$. On the other hand, if $Y$ is an invariant and complemented subspace of $X$, then $Y=\overline{\mathrm{span}}\left(\Gamma_1\right)$ for some $\Gamma_1\subset \Gamma_0$. In order to see that $Y$ is an image of an idempotent invariant operator on $X$ (which has to satisfy $\widehat{T}=\mathbbm{1}_{\Gamma_1}$), we apply the following known lemma due to Rudin.
\begin{lem} \label{rudin}Suppose $X\subset L^1(G)$ is invariant. Then any invariant and complemented subspace $Y\subset X$ is the image of an invariant projection acting on $X$.\end{lem}
\textit{Proof.} Let $P$ be a projection from $X$ onto $Y$ and $Q$ be defined by $Qf=\int_G \left(\tau_{-x}\circ P\circ \tau_x\right)f\d \mu(x)$ in the sense of an $L^1$-valued Riemann integral of a continuous function with respect to the Haar measure $\mu$. By change of variables, $Q\circ \tau_y=\tau_y\circ Q$. For $f\in X$, we have $\tau_xf\in X$, $P\left(\tau_xf\right)\in Y$ and thus $Qf\in Y$. For $f\in Y$, $\tau_xf\in Y$ gives $\left(\tau_{-x}\circ P\circ \tau_x\right)f=\tau_{-x}\left(\tau_xf\right)=f$, so $Q$ is a projection, while the continuity is obvious.\par

It is easy to see that the family of all supports of idempotent multipliers is a Boolean ring (containing all finite sets in the case of compact $G$). We will call it the idempotent ring of $X$. The Cohen-Helson-Rudin theorem states that for $X=L^1(G)$ (in the full locally compact abelian, possibly non-metrizable generality) it is generated by cosets of open subgroups of $\Gamma$. For $f\in H^1(\T)= \overline{\mathrm{span}}\{\e^{\i n t}:n=0,1,2,\ldots\}\subset L^1(\T)$ the idempotent ring is bigger. Indeed, for and any lacunary sequence $(n_k)_{k=1}^\infty$ the Paley inequality \[\left(\sum_{k} \left|\widehat{f}\left(n_k\right)\right|^2\right)^{\frac{1}{2}}\lesssim \left\|f\right\|_1\] holds, with a constant dependent only on $\inf \frac{n_{k+1}}{n_k}$. Thus for any sequence of signs $\varepsilon_1,\varepsilon_2,\ldots$ the map $\sum_n a_n \e^{\i nt}\mapsto \sum_k \varepsilon_k a_{n_k} \e^{\i n_k t}$, which we will denote by $P_{\left(n_k\right),\left(\varepsilon_k\right)}$ is bounded as an $H^1(\T)\rightarrow H^2(\T)$ operator (where $H^2(\T)=H^1(\T)\cap L^2(\T)$ equipped with $L^2$ norm) and consequently as an $H^1(\T)\rightarrow H^1(\T)$ one. This implies that the idempotent ring of $H^1(\T)$  besides arithmetic sequences and finite sets (as restrictions of $L^1(\T)$ multipliers) contains also lacunary sequences. It has been conjectured by Pełczyński and proved by Klemes \cite{Klemes} that they are in fact its generators. \par 
In this section we will be concerned with the space 
\[H^1(\T\times \T)=\overline{\mathrm{span}}\{\e^{\i \left(n_1 t_1+n_2t_2\right)}:n_1,n_2\geq 0\}\subset L^1\left(\T^2\right)\] (not to be confused with $H^1\left(\T^2\right)$, consisting of functions $f\in L^1\left(\T^2\right)$ with coordinatewise Riesz transforms in $L^1\left(\T^2\right)$). Clearly, for any invariant operators  $T_1,T_2$ acting on $H^1(\T)$, there is an invariant operator $T_1\otimes T_2$ acting on $H^1(\T\times \T)$ and satisfying $\widehat{T}\left(n_1,n_2\right)=\widehat{T_1}\left(n_1\right)\widehat{T_2}\left(n_2\right)$. Thus the idempotent ring of $H^1(\T\times \T)$ contains sets of the form $A_1\times A_2$, where $A_1,A_2$ are in the idempotent ring of $H^1(\T)$. Sligthly more general, for any $M\in\mathrm{GL}\left(2,\mathbb{Z}\right)$ the map $S_M: L^1\left(\T^2\right)\rightarrow  L^1\left(\T^2\right)$ defined by $S_M\left(\sum_{n}a_n\e^{\i \langle n,t\rangle}\right)= \sum_n a_n\e^{\i \langle Mn,t\rangle}$ is a (non necessarily surjective) isometry, because $\left(S_Mf\right)(t)=f\left(M^*t\right)$. Under the additional assumption of $M(\mathbb{N}\times\mathbb{N})\subset \mathbb{N}\times\mathbb{N}$, the multiplier of the idempotent operator $S_M^{-1}\circ\left(T_{A_1}\otimes T_{A_2}\right)\circ S_M$ is supported on the set $\left(\mathbb{N}\times\mathbb{N}\right)\cap M^{-1}\left(A_1\times A_2\right)$. 

Another way to produce an idempotent multiplier on $H^1\left(\mathbb{T}\times \mathbb{T}\right)$ out of one dimensional case is to take as $A_1$ and $A_2$ in the above construction the sequences $(2^k)_{k=1}^\infty$ and move any element of the resulting set $A_1\times A_2$ within its dyadic rectangle. One can prove (cf. \cite{Woj}) using the Littlewood - Paley theory, that the sequence obtained in this way is a bounded Fourier multiplier.
Hence the idempotent ring of $H^1\left(\mathbb{T}\times \mathbb{T}\right)$ contains all subsets of $\mathbb{N}\times\mathbb{N}$ whose intersections with dyadic rectangles are of bounded cardinality. The main purpose of this section is to give an example of a Fourier multiplier on $H^1\left(\mathbb{T}\times \mathbb{T}\right)$ that could not be derived by the manipulations described above from one-dimensional results. Namely we prove the following.\par 

\begin{thm}\label{mnoznik} Let $\left(d_k\right)_{k=1}^\infty$, $\left(N_k\right)_{k=1}^\infty$ be sequences of natural numbers such that $\left(d_k\right)_{k=1}^\infty$ is lacunary, $N_k\mid N_{k+1}$, $N_k<N_{k+1}$. Denote \[A_k=\left\{\left(n_1,n_2\right)\in \mathbb{N}\times\mathbb{N}: n_1+n_2=d_k\right\}\] and $B_k=\left\{\left(n_1,n_2\right)\in A_k:N_k\mid n_1\right\}$.  Then the following conditions are equivalent.\\
(i) There exist constants $a,C$ such that $d_{k}\leq CN_{k+a}$.\\
(ii) The set $B=\bigcup_{k=1}^\infty B_k$ is in the idempotent ring of $H^1\left(\mathbb{T}\times \mathbb{T}\right)$. \end{thm}
In order to guarantee that intersections of our $B$ with dyadic rectangles can be arbitrarily large it is enough to assume $\frac{d_k}{N_k}\to\infty $. A canonical example of such sequences is $d_k=N_{k+1}=k!$. \par
\textit{Proof.} The operator $T$ such that $\widehat{Tf}=\mathbbm{1}_B \widehat{f}$ is well defined on polynomials and it suffices to examine its boundedness on them. Let $L=\begin{pmatrix}1& 1\\ 1& -1\end{pmatrix}$.  
Then the isometry
$$
S_L: H^1\left(\mathbb{T}\times\mathbb{T}\right)\to \overline{\mathrm{span}}\left\{ \e^{\i \left(n_1t_1+n_2t_2\right)}:n_1\geq 0, \left|n_2\right|\leq n_1,2\mid n_1+n_2\right\}
$$
identifies $H^1\left(\mathbb{T}\times \mathbb{T}\right)$ with an invariant subspace of 
$$
H^1(\mathbb{T})\otimes L^1(\mathbb{T})=\overline{\mathrm{span}}\left\{ \e^{\i \left(n_1t_1+n_2t_2\right)}:n_1\geq 0\right\}.
$$
The operator $P_{\left(d_k\right),\left(\varepsilon_k\right)}\otimes \mathrm{id}_{L^{1}}$
is a bounded idempotent on this space and therefore $$Q_{\left(\varepsilon_k\right)}= S_L^{-1}\circ \left(P_{\left(d_k\right),\left(\varepsilon_k\right)}\otimes \mathrm{id}_{L^{1}\left(\mathbb{T}\right)}\right)\circ S_L$$ is a bounded idempotent on 
$H^1\left(\mathbb{T}\times \mathbb{T}\right)$. Since $L\left(\bigcup A_k\right)=\left(\left\{d_1,d_2,\ldots\right\}\times \mathbb{Z}\right) \cap L\left(\mathbb{N}\times \mathbb{N}\right)$, the map $Q_{\left(1,1,\ldots\right)}$ acts on $H^1\left(\mathbb{T}\times\mathbb{T}\right)$ as given by $\widehat{Q_{\left(1,1,\ldots\right)}f}=\mathbbm{1}_{\bigcup A_k}\widehat{f}$. 
Since $T=T\circ Q_{(1,1,\dots)}$, it is enough to prove that $T$ is bounded on polynomials whose Fourier transforms are supported in $\bigcup A_k$. \par 

Let $f$ be such polynomial and let $\widehat{f_k}=\mathbbm{1}_{A_k}\widehat{f}$. Then $f=\sum_{k=1}^n f_k$ for some $n$. For any choice of $\varepsilon_1,\varepsilon_2,\ldots\in\{-1,1\}$ we have $Q_{\left(\varepsilon_k\right)}f=\sum_{k=1}^n \varepsilon_k f_k$. Thus $\left\|\sum_{k=1}^n\varepsilon_k f_k\right\|_1\lesssim \|f\|_1$ (since the norm of $P_{\left(d_k\right),\left(\varepsilon_k\right)}$ depends only on $\left(d_k\right)$). 
Since $f=Q_{\left(\varepsilon_k\right)}Q_{\left(\varepsilon_k\right)}f$ we get the reverse estimate and consequently $\left\|\sum_{k=1}^n\varepsilon_k f_k\right\|_1\sim \|f\|_1$. Applying pointwise the Khintchine inequality we get $\|f\|_1\sim \left\|\left(\sum_{k=1}^n\left|f_k\right|^2\right)^{ \frac{1}{2}}\right\|_1$. Moreover, 
$$
\aligned
\left|f_k\right|=&\left|\sum_{n_1+n_2=d_k}\widehat{f_k}\left(n_1,n_2\right)
\e^{\i\left(n_1t_1+n_2t_2\right)}\right|\\
=& \left|\sum_{n_1=0}^{d_k}\widehat{f_k}\left(n_1,d_k-n_1\right)\e^{\i\left(n_1t_1+\left(d_k-n_1\right)t_2\right)}\right|\\
=& \left|\sum_{n_1=0}^{d_k}\widehat{f_k}\left(n_1,d_k-n_1\right)\e^{\i n_1\left(t_1-t_2\right)}\right|.
\endaligned
$$ 
This allows us to write $\|f\|_1\sim \left\|\left(\sum_{k=1}^n\left|\widetilde{f_k}\right|^2\right)^{ \frac{1}{2}}\right\|_1$, where $\widetilde{f_k}\in H^1\left(\mathbb{T}\right)$ is a polynomial of degree at most $d_k$ satisfying $\widehat{\widetilde{f_k}}\left(j\right)=\widehat{f_k}\left(j,d_k-j\right)$. Since $\widehat{Tf_k}\left(j,d_k-j\right)= \widehat{f_k}\left(j,d_k-j\right)$ for $N_k\mid j$ and $0$ otherwise, the action of $T$ on $\widetilde{f_k}$ is given by 
$$
\widetilde{Tf_k}=\sum_{N_k\mid j}\widehat{\widetilde{f_k}}\left(j\right)\e^{\i j t}= \widetilde{f_k}\ast \omega_{N_k},
$$ 
where $\omega_N$ is a measure given by $\widehat{\omega_{N}}\left(j\right)=1$ for $N\mid j$ and $0$ otherwise. Thus $\|Tf\|_1\lesssim\|f\|_1$ transforms into 
$$
\left\|\left(\sum_{k=1}^n\left|\widetilde{f_k}\ast \omega_{N_k}\right|^2\right)^{ \frac{1}{2}}\right\|_1
\lesssim
\left\|\left(\sum_{k=1}^n\left|\widetilde{f_k}\right|^2\right)^{ \frac{1}{2}}\right\|_1.
$$ 
It is easy to see that that the convolution with $\omega_N$ is the orthogonal projection onto the space of $\mathcal{F}^*_{N_k}$-measurable functions, where $\mathcal{F}^*_N$ is the sigma-algebra of subsets of $\mathbb{T}$ whose characteristic functions are $\frac{2\pi}{N}$-periodic. Then the action of $\omega_N$ is given by the formula 
$$
\left(f\ast \omega_N\right)(t)=\frac{1}{N}\sum_{0\leq k<N}f\left(t+\frac{2\pi}{N}k\right).
$$ 
Therefore Theorem \ref{mnoznik} is implied by the following. \par

\begin{thm} \label{stein}Let $d_1,d_2,\ldots$ and $N_1,N_2,\ldots$ be sequences of natural numbers such that $\left(d_k\right)$ is lacunary, $N_k\mid N_{k+1}$ and $N_k<N_{k+1}$. Then the following conditions are equivalent.\\
(i) There exist constants $a,C$ such that $d_{k}\leq CN_{k+a}$.\\
(ii) For any trigonometric polynomials $f_1,\ldots,f_n\in H^1\left(\mathbb{T}\right)$ such that $\deg f_k\leq d_k$ the inequality
\[\E\sqrt{\sum_{k=1}^n\left|f_k\right|^2}\gtrsim \E\sqrt{\sum_{k=1}^n\left|\E^*_{N_k}f_k\right|^2}\]
is satisfied. \\
(iii) For any trigonometric polynomials $f_1,\ldots,f_n\in H^1\left(\mathbb{T}\right)$ such that $\deg f_k\leq d_k$ the inequality
\[\left\|\left(f_k\right)_{k=1}^n\right\|_{L^1\left(\ell^2\right)}\gtrsim \left\|\left(\E^*_{N_k}f_k\right)_{k=1}^n\right\|_{\mathrm{ind}}.\]\end{thm} \par 
The condition $N_k<N_{k+1}$ is here only to reduce the technical difficulties in the proof of implication (ii)$\implies$(i) by guaranteeing the lacunarity of $\left(N_k\right)$. It is easy to see by a standard Khintchine-based linearisation argument that if we allow repetitions in the sequence $\left(N_k\right)$, then the condition (ii) for sequences $\left(d_k\right)$, $\left(N_k\right)$ is equivalent to condition (ii) for sequences $\left(d'_k\right)$, $\left(N'_k\right)$ such that $\left(N'_k\right)$ is $\left(N_k\right)$ with repetitions removed and $d'_k=\max_{N_j=N'_k} d_j$. \par 
\textit{Proof.} We will prove first (i)$\implies$(ii). We do not need $N_k\mid N_{k+1}$. Suppose that (i) is violated and (ii) holds. Take $d'_k=\left\lfloor \frac{d_k}{2}\right\rfloor- \left\lfloor \frac{d_k}{2}\right\rfloor\left(\mathrm{mod}\ N_k\right)$. The condition (i) is violated for sequences $\left(d'_k\right)$ and $\left(N_k\right)$ as well, since $d_k\leq 2\left(d'_k+N_k+1\right)$. The condition (ii) is satisfied for the sequences $\left(d'_k\right)$ and $\left(N_k\right)$ and two sided polynomials $f_k\in L^1(\T)$ (not necessarily in $H^1(\T)$), because for any polynomial $f_k\in L^1(\T)$ of degree $\leq d'_k$ we have $f_k=\e^{-i d'_kt}F_k$, where $F_k$ is a polynomial in $H^1(\T)$ of degree $\leq d_k$ and $\E^*_{N_k}$ commutes with multiplying by $\e^{-id'_k t}$. \par 
Let us take an arbitrary $a$. Then there exists $k$ such that $d'_k>N_{k+a}$. We may apply the inequality from condition (ii) to $f_k=\ldots=f_{k+a}=K_{d'_k}$ and $f_j=0$ for $j\notin\{k,\ldots,k+a\}$, where $K_{d'_k}$ is the Fej\'{e}r kernel. Let $I_k=\left[-\frac{\pi}{d'_k},\frac{\pi}{d'_k}\right]$. It is easy to check that $K_{d'_k}\gtrsim d'_k \mathbbm{1}_{I_k}$ and that the condition $d'_k>N_{k+a}$ implies that for $k\leq j\leq k+a$ the function $\E^*_{N_j}\mathbbm{1}_{I_k}$ is supported on a sum of $N_j$ disjoint intervals of length $\frac{2\pi}{d'_k}$ and attains only vlaues $0$ and $\frac{1}{N_j}$. We have \[(a+1)^\frac12=\E\sqrt{\sum_{j=k}^{k+a}\left|K_{d'_k}\right|^2}\gtrsim \E\sqrt{\sum_{j=k}^{k+a}\left|\E^*_{N_j}K_{d'_k}\right|^2}\gtrsim d'_k\E\sqrt{\sum_{j=k}^{k+a}\left|\E^*_{N_j}\mathbbm{1}_{I_k}\right|^2}\]\[\geq d'_k\sum_{j=k+1}^{k+a}\frac{1}{N_j}\mu\left(\mathrm{supp}\ \E^*_{N_j}\mathbbm{1}_{I_k}\setminus \mathrm{supp}\ \E^*_{N_{j-1}}\mathbbm{1}_{I_k}\right)\geq \]\[d'_k\sum_{j=k+1}^{k+a}\frac{1}{N_j}\left(\frac{2\pi N_j}{d'_k}-\frac{2\pi N_{j-1}}{d'_k}\right)\gtrsim a,\] where the last inequality follows from lacunarity of $\left(N_k\right)$. Since $a$ was arbitrary, this is a contradiction completing the proof of implication (ii)$\implies$(i).\par 
Since Corollary \ref{2indgeqL1l2} gives (iii)$\implies$(ii) immediately, the only thing left is the proof of implication (i)$\implies$(iii). We can assume without loss of generality that $d_k\leq CN_{k+1}$, because in general case we consider sequences $\left(d_{ak+b}\right)_{k=1}^\infty$, $\left(N_{ak+b}\right)_{k=1}^\infty$ for $b=0,\ldots,a-1$ and add the resulting inequalities. \par 
We will need two lemmas. We will assume that all functions involved are absolutely continuous. \par 

\begin{lem} \label{conv}
Let $\psi_N(x)=x-\frac{2\pi}{N}k$ whenever $\frac{2\pi}{N}k\leq x<\frac{2\pi}{N}(k+1)$. Then $\left(\E-\E^*_N\right)f=f'\ast\psi_N$. 
\end{lem}

\textit{Proof.} Let $x_k= x+\frac{2\pi}{N}k$. Then 
$$
\aligned
\left(\E f-\E^*_N f\right)(x) &= \frac{1}{2\pi}\sum_{0\leq k<N}\int_{x_k}^{x_{k+1}}\left(f(x)-f\left(x_k\right)\right)\d x \\
&= \frac{1}{2\pi}\sum_{0\leq k<N}\int_{x_k}^{x_{k+1}}\int_{x_k}^x f'(y)\d y\d x \\
&= \frac{1}{2\pi}\sum_{0\leq k<N}\int_{x_k}^{x_{k+1}}\left(x_{k+1}-y\right)f'(y)\d y \\
&= \frac{1}{2\pi}\int \psi_N(x-y)f'(y)\d y=\left(\psi_N\ast f'\right)(x).
\endaligned
$$\par

\begin{deff}
For any function $f:\mathbb{T}\to\mathbb{C}$ we put
$$
\Ber f=\frac{\E \left|f'\right|}{\E|f|}.
$$
\end{deff}

\begin{lem}\label{ber} Let $\left\|\cdot\right\|$ be a differentiable norm on $\mathbb{C}^n$ such that $\left\|\left(x_1,\ldots,x_n\right)\right\|\leq\left\|\left(x_1,\ldots,x_n\right)\right\|_{\ell^1_n}$ and $\left\|\left(x_1,\ldots,x_n\right)\right\|=\left\|\left(\left|x_1\right|,\ldots,\left|x_n\right|\right)\right\|$ for any $x_1,\ldots,x_n\in\mathbb{C}$. Then
\[\Ber\left\|\left(f_1,\ldots,f_n\right)\right\|\leq \left\|\left(\Ber f_1,\ldots,\Ber f_n\right)\right\|_{\ast},\]
where $\left\|v\right\|_{\ast}=\sup_{\|w\|=1}\left|\langle v,w\rangle\right|$ is the dual norm. 
\end{lem}

\textit{Proof.} Denote $f=\left(f_1,\ldots,f_n\right)$. Then by the pointwise inequality $\left|\|f\|'\right|\leq \|f'\|$,
$$
\aligned
\Ber \|f\|&=\frac{\E \left|\|f\|'\right|}{\E\|f\|} 
\leq \frac{\E \left\|f'\right\|}{\E \|f\|}
\leq \frac{\E \left\|f'\right\|_1}{\E \|f\|}\\
&= \frac{\sum_k \E\left|f_k\right|\Ber f_k}{\E\|f\|}\\
&\leq\frac{\left\|\left(\E\left|f_1\right|,\ldots,\E\left|f_n\right|\right)\right\| \cdot \left\|\left(\Ber f_1,\ldots,\Ber f_n\right)\right\|_{\ast}}{\E\|f\|}\\
&\leq \left\|\left(\Ber f_1,\ldots,\Ber f_n\right)\right\|_{\ast}.
\endaligned
$$\par 

In order to prove (i)$\implies$(ii), we follow the usual routine of Section 1. The divisibility condition $N_k\mid N_{k+1}$ implies that $\left(\mathcal{F}^*_{N_k}\right)_{k=1}^n$ is a decreasing filtration. Denote $\varphi_k=\E^*_{N_k}f_k$ and $\lambda_0=0$, $\lambda_{k}=\E^*_{N_{k+1}}\sqrt{\left|\varphi_k\right|^2+\lambda_{k-1}^2}$ for $k\geq 1$. Obviously $\Ber \E^*_{k+1}f\leq \Ber f$ for nonnegative $f$. Note that $\varphi_k$ is a polynomial of degree at most $d_k$. By Lemma \ref{ber} for $\|\cdot\|_{\ell^2_2}$ and Bernstein inequality,
$$
\Ber \lambda_k\leq \sqrt{\left(\Ber \varphi_k\right)^2+ \left(\Ber\lambda_{k-1}\right)^2}\leq \sqrt{d_{k}^2+\left(\Ber\lambda_{k-1}\right)^2}
$$
 and thus, by lacunarity, $\Ber\lambda_k\leq \sqrt{d_k^2+\ldots+d_1^2}\lesssim d_k$. Therefore, by Lemma \ref{conv},
$$
\aligned
\lambda_k&=\E^*_{N_{k+1}}\lambda_k\\
&\leq\E\lambda_k+\left|\E\lambda_k-\E^*_{N_{k+1}}\lambda_k\right|\\
&= \E\lambda_k + \left|\psi_{N_{k+1}}\ast \lambda_k'\right|\\
&\leq\E\lambda_k +\frac{2\pi}{N_{k+1}}\E\left|\lambda_{k}'\right|\\
&=\left(1+\frac{2\pi}{N_{k+1}}\Ber\lambda_k\right)\E\lambda_k\\
&\lesssim \left(1+\frac{2\pi d_k}{N_{k+1}}\right)\E\lambda_k\\
&\lesssim\E\lambda_k.
\endaligned
$$
Ultimately, applying Lemma \ref{telescoping} to $\varphi_k$ and Lemma \ref{inf} to $f_k$ and $\E^*_{N_k}$ we get \[ \left\|\left(\varphi_k\right)_{k=1}^n\right\|_{\mathrm{ind}}\lesssim \sqrt{\E\lambda_n}\sqrt{\max_{k\leq n}\sup\lambda_k}\lesssim \sqrt{\E\lambda_n}\sqrt{\max_{k\leq n}\E\lambda_k}\]
\[=\E\lambda_n \leq\E\sqrt{\sum_{k=1}^n \left|f_k\right|^2}\] completing the proof of the theorem.\par 
It has to be noted that all the information about $f_1,\ldots,f_n$ we used in the proof of the crucial implication (i)$\implies$(iii) was contained in $\Ber \varphi_1,\ldots,\Ber \varphi_n$. Thus in the (i)$\implies$(ii) implications in the Theorem \ref{stein} and Corollary \ref{renorm} we only need $\Ber \varphi_k\lesssim N_{k+a}$, without the assumption that $\varphi_k$ are polynomials or belong to $H^1(\T)$.\par 
The main part of the following corollary is somewhat trivial. We state it separately for its similarity with Theorem \ref{l1renorm} and Theorem \ref{transference} (compare results \cite{meyerdechamp}, \cite{krystian}) and because we find the fact that the assumption $d_k\lesssim N_{k+a}$ can be weakened in neither Theorem \ref{stein} nor Corollary \ref{renorm} worth a proof. 
\begin{cor}\label{renorm}Let $d_1,d_2,\ldots$ and $N_1,N_2,\ldots$ be sequences of natural numbers such that $\left(d_k\right)$ is lacunary and $N_k\mid N_{k+1}$. Then the following conditions are equivalent.\\
(i) There exist constants $a,C$ such that $d_{k}\leq CN_{k+a}$.\\
(ii) For any polynomials $\varphi_1,\ldots,\varphi_n\in H^1(\T)$ such that $\varphi_k$ is $\mathcal{F}^*_{N_k}$-measurable and $\deg \varphi_k\leq d_k$, the norms $\left\|\left(\varphi_k\right)_{k=1}^n\right\|_{L^1\left(\ell^2\right)}$ and $\left\|\left(\varphi_k\right)_{k=1}^n\right\|_{\mathrm{ind}}$ are comparable.\end{cor}
\textit{Proof.} The implication (i)$\implies$(ii) is almost trivial and does not need $\varphi_k\in H^1(\T)$: the inequality $\left\|\left(\varphi_k\right)_{k=1}^n\right\|_{L^1\left(\ell^2\right)}\lesssim\left\|\left(\varphi_k\right)_{k=1}^n\right\|_{\mathrm{ind}}$ is just Corollary \ref{2indgeqL1l2} and the inequality $\left\|\left(\varphi_k\right)_{k=1}^n\right\|_{L^1\left(\ell^2\right)}\gtrsim\left\|\left(\varphi_k\right)_{k=1}^n\right\|_{\mathrm{ind}}$ is just Theorem \ref{stein} applied to $f_k=\varphi_k$.\par In order to prove the implication (ii)$\implies$(i) we proceed just as in the proof of Theorem \ref{stein}. We assume for the sake of contradiction that (i) is violated and (ii) is true and then modify sequence $\left(d_k\right)$ in such a way that (ii) is satisfied without the restriction $\varphi_k\in H^1(\T)$. Then we choose arbitrary $a,C\in\mathbb{N}$ and find $k$ such that $d_k> CN_{k+a}$. Let us put $\varphi_{k}=\ldots=\varphi_{k+a}=\E^*_{N_{k+a}}K_{CN_{k+a}}$ and $\varphi_j=0$ for $j\notin\{k,\ldots,k+a\}$. It is easy to see that $\E^*_{N_{k+a}}K_{CN_{k+a}}$ has the same distrubution as $K_C$. Applying (ii) gives \[1=(a+1)^{-\frac12}\left\|\left(\varphi_k\right)_{k=1}^n\right\|_{L^1\left(\ell^2\right)}\gtrsim (a+1)^{-\frac12}\left\|\left(\varphi_k\right)_{k=1}^n\right\|_{\mathrm{ind}}=\E\sqrt{(a+1)^{-\frac12}\sum_{j=k}^{k+a}\varphi_j^2}.\] We replace the function $t\mapsto t^\frac{1}{2}$ by $t\mapsto \min\left(t^\frac12,M\right)$, take the limit with $a\to\infty$ by the law of large numbers and pass with $M$ to infinity by dominated convergence to get $1\gtrsim \sqrt{\E K_C^2}$, which is false when $C$ is large enough, giving the desired contradiction. \par 
We may now combine the results of this section. 

\begin{cor}The independent sum $\left(\bigoplus_{k=1}^\infty H^1_k(\T)\right)_{\mathrm{ind}}$ is isomorphic to a complemented and invariant subspace of $H^1(\T\times\T)$. \end{cor}
\textit{Proof.} Take for instance $d_k=k!=N_{k+1}$. Using the notation of Theorem \ref{mnoznik}, the space $X=\overline{\mathrm{span}}\left\{e^{i\langle n,t\rangle}:n\in B\right\}$ is invariant and complemented in $H^1(\T\times\T)$. The desired isomorphism $T:\left(\bigoplus H^1_k(\T)\right)_{\mathrm{ind}}\to X$ is given by the formula \[T\left(\left(f_k\right)_{k=1}^\infty\right)= \sum_k \sum_{\substack{0\leq j\leq d_k\\N_k\mid j}}\widehat{f_k}\left(\frac{j}{N_k}\right)\e^{i\left\langle\left(j,d_k-j\right),t\right\rangle}\]
where only finitely many $f_k$ are nonzero. Clearly $T$ is one-to-one and has a dense image. Moreover, by the proof of Theorem \ref{mnoznik} and Corollary \ref{renorm}, \[\left\|T\left(\left(f_k\right)_{k=1}^\infty\right)\right\|_{H^1(\T\times\T)}\sim \left\|\left(\sum_{\substack{0\leq j\leq d_k\\N_k\mid j}}\widehat{f_k}\left(\frac{j}{N_k}\right)\e^{ijt}\right)_{k=1}^\infty\right\|_{L^1\left(\ell^2\right)}=\]\[\left\|\left(f_k\left(N_k t\right)\right)_{k=1}^\infty\right\|_{L^1\left(\ell^2\right)}\sim \left\|\left(f_k\left(N_k t\right)\right)_{k=1}^\infty\right\|_{\mathrm{ind}}=\indnorm{\left(f_k\right)_{k=1}^\infty}\] since $f_k(t)$ and $f_k\left(N_k t\right)$ have the same distribution, which completes the proof. \par 
Although the assumption $N_k\mid N_{k+1}$ is crucial for our proof of Theorem \ref{stein}, it turns out that at the price of a stronger qualitative assumption on the sequences $\left(d_k\right)$ and $\left(N_k\right)$ we can retain the comparability of norms proved in Corollary \ref{renorm}. 

\begin{thm}\label{l1renorm}Let $d_1,d_2,\ldots$ and $N_1,N_2,\ldots$ be sequences of natural numbers such that $\left(d_k\right)$ is lacunary and $\sum_k\frac{d_k}{N_{k+a}}<\infty$ for some $a$. Then for any polynomials $\varphi_1,\ldots,\varphi_n\in H^1(\T)$ such that $\varphi_k$ is $\mathcal{F}^*_{N_k}$-measurable (we do not imply that $\left(\mathcal{F}^*_{N_k}\right)_k$ is a filtration) and $\deg \varphi_k\leq d_k$, the norms $\left\|\left(\varphi_k\right)_{k=1}^n\right\|_{L^1\left(\ell^2\right)}$ and $\left\|\left(\varphi_k\right)_{k=1}^n\right\|_{\mathrm{ind}}$ are comparable.\end{thm} 
\textit{Proof.} Without loss of generality we may assume $a=1$ and $2\pi\sum \frac{d_k}{N_{k+1}}<\frac12$. Let us fix $k$ and $t_{k+1},\ldots,t_n$ and denote $u_k(t)=\sqrt{\sum_{j=1}^{k-1}\left|\varphi_j(t)\right|^2+ \sum_{j=k+1}^n\left|\varphi_j(t_j)\right|^2}$. By Fej\'{e}r-Riesz lemma, $u_k$ is the modulus of a trigonometric polynomial of degree $\leq d_{k-1}$. It follows from Lemma \ref{conv} that \[\left|\int_{\T}\sqrt{u_k^2(t)+\left|\varphi_k(t)\right|^2}\frac{\d t}{2\pi}-\int_{\T\times \T}\sqrt{u_k^2(s)+\left|\varphi_k(t)\right|^2}\frac{\d s}{2\pi}\frac{\d t}{2\pi}\right|=\]
\[\left|\int_{\T}\frac{1}{N_k}\sum_{j=0}^{n-1}\sqrt{u_k^2\left(t+\frac{2\pi}{N_k}j\right)+\left|\varphi_k(t)\right|^2}\frac{\d t}{2\pi}-\int_{\T\times \T}\sqrt{u_k^2(s)+\left|\varphi_k(t)\right|^2}\frac{\d s}{2\pi}\frac{\d t}{2\pi}\right|\]
\[=\left|\int_{\T}\left(\E^*_{N_k}\sqrt{u_k^2+\left|\varphi_k(t)\right|^2}\right)(t)\frac{\d t}{2\pi}-\int_{\T}\E\sqrt{u_k^2+\left|\varphi_k(t)\right|^2}\frac{\d t}{2\pi}\right|\]\[\leq \int_{\T}\left|\left(\E-\E^*_{N_k}\right)\sqrt{u_k^2+\left|\varphi_k(t)\right|^2}\right|(t)\frac{\d t}{2\pi}\leq \]\[\int_{\T}\frac{2\pi}{N_k}\E\left|\left(\sqrt{u_k^2+\left|\varphi_k(t)\right|^2}\right)'\right|\frac{\d t}{2\pi}\leq \int_{\T}\frac{2\pi}{N_k}\E\left|u_k'\right|\frac{\d t}{2\pi}=\]\[\frac{2\pi\Ber u_k}{N_k}\E u_k\leq \frac{2\pi d_{k-1}}{N_k}\E u_k.\]
By the above inequality and Corollary \ref{2indgeqL1l2}, \[\left|\Lnorm{\left(\varphi_k\right)_{k=1}^n}-\indnorm{\left(\varphi_k\right)_{k=1}^n}\right|= \] \[\left|\sum_{k=1}^n \int_{\T^{n-k}}\left( \int_{\T}\sqrt{\left|\varphi_k(t)\right|^2+u_k^2(t)}\frac{\d t}{2\pi} +\right.\right.\]\[\left.\left.- \int_{\T\times \T}\sqrt{\left|\varphi_k(t_k)\right|^2+u_k^2(t)}\frac{\d t}{2\pi}\frac{\d t_k}{2\pi}\right)\frac{\d t_{k+1}}{2\pi}\ldots \frac{\d t_{k+1}}{2\pi}\right|\]
\[\leq \sum_{k=1}^n \frac{2\pi d_{k-1}}{N_k}\E u_k\leq 2\pi\sum \frac{d_{k}}{N_{k+1}}\indnorm{\left(\varphi_k\right)_{k=1}^\infty}<\frac12\indnorm{\left(\varphi_k\right)_{k=1}^\infty}\]
which gives the desired comparability.

\section{Independent sum of $H^1_n(\delta)$}

Let us fix the notation concering the Cantor group and briefly recall the basic information about the dyadic Hardy space. A standard reference in this matter is \cite{bigmuller}. \par 
Let $\mathcal{F}_n$ be the sigma-algebra generated by dyadic subintervals of $[0,1]$ of length $2^{-n}$ and $\mathcal{F}^*_n$ be the sigma-algebra of subsets of $[0,1]$ whose characteristic functions are $2^{-n}$-periodic. It is easy to see that $\mathcal{F}_n$ and $\mathcal{F}^*_n$ are independent. \par 
We will work with the Cantor group $\mathbb{Z}_2^{\mathbb{N}^+}$, where our model for $\mathbb{Z}_2$ is $\{0,1\}$. Its elements can be identified with subsets of $\mathbb{N}^+$. Elements of the dual (Walsh) group $\widehat{\mathbb{Z}_2^{\mathbb{N}^+}}$ are finitely supported $\mathbb{Z}_2$-valued sequences indexed by $\mathbb{N}^+$ and can be identified with finite subsets of $\mathbb{N}^+$. Thus $\max A$ and $\min A$ are well defined for $A$ in the Walsh group and $\min x$ is well defined for $x$ in the Cantor group (with the natural convention $\max \emptyset=0, \min\emptyset=\infty$). The duality is given by the formula $w_A(x)=(-1)^{\sum_{k=1}^\infty a_k x_k}$. Addition on any of these groups corresponds to the symmetric difference of sets. The metric on the Cantor group will be $d\left(x,y\right)=d\left(x-y,\emptyset\right)=2^{-\min\left(x-y\right)}=2^{-\min\left\{k:x_k\neq y_k\right\}}$.\par 
We will frequently use the canonical measure-preserving identification of $[0,1]$ and the Cantor group given by the formula \[\mathbb{Z}_2^{\mathbb{N}^+}\ni x\mapsto \sum_{k=1}^\infty 2^{-k}x_k\in [0,1].\] In this setting, $w_A=\prod_{k\in A}r_k$, where $r_k(t)=\mathrm{sign}\sin\left(2^k\pi t\right)$ are the Rademacher functions. The Walsh functions form a complete orthonormal system in $L^2[0,1]$ just as they did in $L^2\left(\mathbb{Z}_2^{\mathbb{N}^+}\right)$. Also, $\mathcal{F}_n$ corresponds to the product of full sigma-algebra of subsets of $\mathbb{Z}_2^{\{1,\ldots,n\}}$ and the trivial sigma-algebra of subsets of $\mathbb{Z}_2^{\{n+1,\ldots\}}$, while $\mathcal{F}^*_n$ is the same with the words "trivial" and "full" interchanged. The metric on the Cantor group corresponds, up to scaling, to the dyadic metric on $[0,1]$. \par
We denote the standard martinagle differences by $\Delta_0=\E_0$, $\Delta_k=\E_k-\E_{k-1}$ for $k\geq 1$ and define the square function by \[S(f)=\left(\sum_{n=0}^\infty\left|\Delta_n f\right|^2\right)^{\frac12}\] An equivalent definition is \[S(f)=\left(\sum_{n=0}^\infty\left|\sum_{\max A=n}\left\langle f,w_A\right\rangle w_A\right|^2\right)^{\frac{1}{2}},\] since the space of $\mathcal{F}_n$-measurable functions is $\mathrm{span}\left\{w_A:A\subset\{1,\ldots,n\}\right\}$ and thus the space of $n$-th martingale diffrences is exactly $\mathrm{span}\left\{w_A:\max A=n\right\}$. The dyadic Hardy space $H^1\left(\delta\right)$ is the space of integrable functions defined on $[0,1]$ such that the norm $\|f\|_{H^1(\delta)}=\E S(f)$ is finite, equipped with this norm. It is known to be a Banach space, as a consequence of the Davis inequality and Doob martingale convergence theorem. By $H^1_n(\delta)$ we mean the subspace of $H^1\left(\delta\right)$ consisting of $\mathcal{F}_n$-measurable functions. The first definition of the square function gives a canonical isometry $\iota: H^1\left(\delta\right)\to L^1\left([0,1],\ell^2\right)$ by $\iota(f)=\left(\Delta_n f\right)_{n=0}^\infty$. The range of $\iota$ is exactly the space $\left\{\left(f_n\right)\in L^1\left(\ell^2\right): f_n=\Delta_nf_n\right\}$ (since the latter is a Banach space as well). It is convenient to consider an orthonormal basis in $L^2\left(\ell^2\right)$ consisting of functions of the form $w_A\otimes e_n$, where $e_n$ are unit vectors in $\ell^2$. Then $\iota\left(H^1(\delta)\right)$ is the closure of the span of $w_A\otimes e_{\max A}$, where $A$ runs through finite subsets of $\mathbb{N}^+$. We will denote by $H^2(\delta)$ the function space on $[0,1]$ with the norm $\|f\|_{H^2(\delta)}=\left(\E\left(Sf\right)^2\right)^\frac12= \left\|\iota(f)\right\|_{L^2\left(\ell^2\right)}$. An easy calculation shows that $\iota$ is an $\|\cdot\|_{L^2}\to \|\cdot\|_{L^2\left(\ell^2\right)}$ isometry, so $H^2(\delta)=\cap L^2$. However, we will sometimes use the $H^2(\delta)$ notation to emphasize that we are dealing with the norm expressed in terms of the square function.\par 
\par 
The following theorem of M\"{u}ller and Schechtman \cite{bigmuller}, \cite{smallmuller} introduces the space and $\left(\sum_{n=1}^\infty H^1_n(\delta)\right)_{\mathrm{ind}}$. We present a new proof which gives a better constant.\par 
\begin{thm}\label{indcompl} Let $\left(I_n\right)_{n=1}^\infty$ be a family of intervals in $\mathbb{N^+}$ such that $\max I_n <\min I_{n+1}$ and let $G\left(I_n\right)$ be $\mathrm{span}\left\{w_A:A\subset I_n\right\}$. Then the orthogonal projection $P$ onto the subspace $\overline{\mathrm{span }}\bigcup G\left(I_n\right)$ extends to an bounded $H^1(\delta)\to H^1(\delta)$ operator.\end{thm}
\textit{Proof.} 
It is enough to prove boundedness of $P$ on functions $f$ which are finite linear combinations of Walsh functions. Denote $I_n=\left[a_n,b_n\right]$.\par 
\textit{Claim.} $\Delta_k Pf=\E^*_{a_n-1}\Delta_k f$ if $k\in I_n$ and $\Delta_k Pf=0$ if $k\notin \bigcup_n I_n$. \par 
Indeed, suppose first that $k\notin \bigcup_n I_n$ and take $A$ such that $A\subset I_n$ for some $n$. Since $\Delta_k$ is the projection onto $\mathrm{span}\{w_A:\max A=k\}$ and $\max A\in I_n$, we have $\Delta_k w_A=0$, so $\Delta_k$ restricted to the image of $Pf$ is identically $0$. On the other hand, if $k\in I_{n_0} $ for some $n_0$. By the identification of $[0,1]$ with $\mathbb{Z}_2^{\mathbb{N}^+}$, the operator $\E^*_{a_n-1}$ is the projection onto $\overline{\mathrm{span}}\left\{w_A:\min A\geq a_{n}\right\}$. Thus, operators $\Delta_kP$ and $\E^*_{a_n-1}\Delta_k$ both correspond to idempotent Fourier multipliers on the Walsh group. We have $\widehat{\Delta_k P}(A)= 1$ iff $\max A=k$ and $A\subset I_n$ for some $n$, and $\widehat{\E^*_{a_{n_0}-1} \Delta_k}(A)=1$ iff $\max A=k$ and $\min A\geq a_{n_0}$ and since $n_0$ is the only $n$ such that $k\in I_{n}$, these conditions on $A$ are equivalent, which proves the claim. \par 
Let $g_n=\left(\sum_{k\in I_n}\left|\Delta_kf_k\right|^2\right)^\frac{1}{2}$. Then $g_n$ is $\mathcal{F}_{b_n}\subset\mathcal{F}_{a_{n+1}-1}$-measurable. Thus, it is independent on $\mathcal{F}^*_{a_{n+1}-1}$ and by Corollary \ref{duallep} applied for $g_n$ and $\mathcal{F}^*_{a_n-1}$, we obtain 
\[\|f\|_{H^1(\delta)}\geq\E\left(\sum_n \sum_{k\in I_n}\left|\Delta_k f_k\right|^2\right)^\frac{1}{2}= \E\left(\sum_n g_n^2\right)^\frac{1}{2}\geq \frac{1}{2}\E\left(\sum_n\left|\E^*_{a_{n}-1}g_n\right|^2\right)^\frac{1}{2}
\geq\] \[\geq \frac{1}{2}\E\left(\sum_n \sum_{k\in I_n} \left|\E^*_{a_n-1}\Delta_k f\right|^2\right)^\frac{1}{2}= \frac{1}{2}\E\left(\sum_n\sum_{k\in I_n} \left|\Delta_k Pf\right|^2\right)^\frac{1}{2}=\frac{1}{2}\left\|Pf\right\|_{H^1(\delta)}.\]
In order to exibit the connection between the range of the projection in Theorem \ref{indcompl} we need the following observation.
\begin{lem}\label{trywdil}Consider the dilation operator $Tf(x)=f\left(2x\left(\mathrm{mod}\ 1\right)\right)$. Then $Tw_A=w_{A+1}$ and $T(fg)=(Tf)(Tg)$. Moreover, $T$ preserves distribution and $STf=TSf$. In particular, $T$ is an isometry in the $\|\cdot\|_{H^1(\delta)}$ norm.\end{lem}
\textit{Proof.} The property $T(fg)=(Tf)(Tg)$ and preserving distribution are obvious. In the Cantor group seting, $T$ corresponds to the superposition with the map $\left(x_1,x_2,\ldots\right)\mapsto \left(x_2,x_3,\ldots\right)$ and thus $Tr_k=r_{k+1}$ and consequently $Tw_A=w_{A+1}$. Therefore \[\left(STf\right)^2=\sum_{n}\left|\sum_{\max A=n}\left\langle Tf,w_A\right\rangle w_A\right|^2= \sum_n\left|\sum_{\max A=n}\left\langle Tf,w_{A+1}\right\rangle w_{A+1}\right|^2= \]\[\sum_n\left|\sum_{\max A=n}\left\langle f,w_{A}\right\rangle w_{A+1}\right|^2= T\left(Sf\right)^2.\]
\par 
Suppose that $\left|I_n\right|=n$ in the setting of Theorem \ref{indcompl}. Let $T_n=T^{a_n-1}$ be the power of $T$ sending $H^1_n(\delta)$ to $G\left(I_n\right)$. For $n\neq m$ the elements of $G\left(C_n\right)$ and $G\left(C_m\right)$ are independent as random variables, since they are built of disjoint sets of Rademacher functions. For these reasons, the space $\left(\mathrm{span }\bigcup G\left(C_n\right),\|\cdot\|_{H^1(\delta)}\right)$ is usually denoted by $\left(\sum_{n=1}^\infty H^1(\delta)_n\right)_{\mathrm{ind}}$. The notational difference highlights the fact that it is not precisely an independent sum in our sense. However, we have the following
\begin{cor}The space $\left(\bigoplus_{n=1}^\infty\iota\left(H^1_n(\delta)\right)\right)_{\mathrm{ind}}$ is isometric to $\left(\sum_{n=1}^\infty H^1_n(\delta)\right)_{\mathrm{ind}}$.\end{cor}
\textit{Proof.} Note that for $n\neq m$, the square functions $S\left(T_n f_n\right)$, $S\left(T_m f_m\right)$ are independent and $S\left(\sum T_n f_n\right)^2=\sum S\left(T_n f_n\right)^2$, both due to the fact that $T_n f_n$ and $T_m f_m$ use disjoint sets of Rademacher functions. By Lemma \ref{trywdil}
\[\aligned \left\|\sum T_n f_n\right\|_{H^1(\delta)}&= \E S\left(\sum T_n f_n\right)\\ &=\E \left(\sum S\left(T_n f_n\right)^2\right)^{\frac12}\\ &= \int_{[0,1]^\infty} \left(\sum S\left(T_n f_n\right)\left(x_n\right)^2\right)^{\frac{1}{2}}\d x_1 \d x_2 \ldots\\ &= \int_{[0,1]^\infty} \left(\sum S\left(f_n\right)\left(x_n\right)^2\right)^{\frac{1}{2}}\d x_1 \d x_2 \ldots\\ &=\int_{[0,1]^\infty} \left(\sum \left\|\iota\left(f_n\right)\left(x_n\right)\right\|_{\ell^2}^2\right)^{\frac{1}{2}}\d x_1 \d x_2 \ldots\\ &=\left\|\left(\iota\left(f_n\right)\right)_{n=1}^\infty\right\|_{\mathrm{ind}}\endaligned\]
which gives an isometry between $\left(\mathrm{span}\bigcup G\left(C_n\right),\|\cdot\|_{H^1(\delta)}\right)$ and $\left(\bigoplus \iota\left(H^1_n(\delta)\right)\right)_{\mathrm{ind}}$. \par

\section{The $K$-closedness result}
We recall certain notions from interpolation theory due to Peetre and Pisier \cite{peetre}, \cite{pisier}.
\begin{deff}Suppose that Banach spaces $X_0,X_1$ are embedded in a linear topological space $X$ (we will call them then a compatibile couple). For $t>0$ and $f$ in the algebraic sum $X_0+X_1$ we define the $K$-functional as \[K_t\left(f,X_0,X_1\right)=\inf_{\substack{f=g+h\\g\in X_0,h\in X_1}}\|g\|_{X_0}+\|h\|_{X_1}.\]\end{deff}

\begin{deff}Let $\left(X_0,X_1\right)$ be a compatible couple and $Y_i$ be a subspace of $X_i$ for $i=0,1$ (we will call $(\left(Y_0,Y_1\right)$ a subcouple of $\left(X_0,X_1\right)$). We will say that the couple $\left(Y_0,Y_1\right)$ is $K$-closed in the couple $\left(X_0,X_1\right)$ if \[K_t\left(f,Y_0,Y_1\right)\leq C K_t\left(f,X_0,X_1\right),\]
where $C$ does not depend on $t$ and $f$. \end{deff}
\begin{lem}\label{kcleq}Let $\left(Y_0,Y_1\right)$ be a subcouple of $\left(X_0,X_1\right)$ and denote the norms on $X_i$ by $\|\cdot\|_i$. The following are equivalent.\\
(i) $\left(Y_0,Y_1\right)$ is $K$-closed in $\left(X_0,X_1\right)$.\\
(ii) For any $f\in Y_0+Y_1$ and a decomposition $f=x_0+x_1$, where $x_i\in X_i$, there exists a decomposition $f=y_0+y_1$ such that $y_i\in Y_i$ and $\|y_i\|_i\leq C\|x_i\|_i$.\end{lem}
\textit{Proof.} The implication (ii)$\implies$(i) is obvious. We will prove (i)$\implies$(ii). Let $f=x_0+x_1$. If $\|x_0\|_0=0$ or $\|x_1\|_1=0$, there is nothing to prove. Otherwise let $t=\frac{\|x_0\|_0}{\|x_1\|_1}$. By $K$-closedness, there exist $y_0,y_1$ such that $f=y_0+y_1$ and $\|y_0\|_0+t\|y_1\|_1\leq 2K_t\left(f,Y_0,Y_1\right)\leq 2C 2K_t\left(f,X_0,X_1\right)\leq 2C\left(\|x_0\|_0+t\|x_1\|_1\right)$ and the result follows. \par 
Let us recall the most elementary version of vector-valued Calder\'{o}n-Zygmund theorem, in which we assume $L^2$-boundedness and existence of the principal value. 
\begin{thm}\label{czo}Let $B_1$, $B_2$ be a reflexive and separable Banach spaces and $\left(X,d,\mu\right)$ be a space of homogeneous type. Let $K$ be a function on $X\times X$ with values in the space $\mathcal{L}(B_1,B_2)$ of bounded operators from $B_1$ to $B_2$, bounded outside any neighbourhood of the diagonal. Suppose that $K$ satifies the H\"{o}rmander condition \[\sup_{r>0}\sup_{d\left(y,y_0\right)<r}\int_{X\setminus B\left(y_0,2r\right)}\left\|K(x,y)-K\left(x,y_0\right)\right\|_{\mathcal{L}\left(B_1,B_2\right)}\d \mu (x)<C_1<\infty.\]
If the operator $T$ acting on $L^1\left(X,B_1\right)$ is given by the formula \[Tf(x)=\pv \int_X K(x,y)f(y)\d \mu(y)\] for any $x\in X$ and $\left\|T\right\|_{L^2\left(B_1\right)\to L^2\left(B_2\right)}<C_2$, then $T$ is of weak type $1-1$ and consequently is bounded from $L^p\left(B_1\right)$ to $L^p\left(B_2\right)$ for $1<p<\infty$. Moreover, $\left\|Tf\right\|_{L^{1,\infty}\left(B_2\right)}\leq C \|f\|_{L^1\left(B_1\right)}$, where $C$ depends only on $C_1,C_2$ and the space $\left(X,d,\mu\right)$.\end{thm}
Our main tool will be the following result of Bourgain \cite{bourginterp}. 
\begin{lem}\label{bourgkcl}\emph{(Bourgain)} Let $\left(X,d,\mu\right)$ be a space of homogeneous type and $S$ be some linear space of $\ell^2$-valued bounded measurable functions on $X$. Denote $S_p=\overline{S}^{\|\cdot\|_{L^p\left(\ell^2\right)}}$. Suppose that the orthogonal projection $P:L^2\to S_2$ is given off-diagonal by a $\mathcal{L}\left(\ell^2\right)$-valued Calder\'{o}n-Zygmund kernel $K$, i.e. \[Pf(x)=\int_{X}K(x,y)f(y)\d \mu(y),\]
for $x\notin\mathrm{supp}\ f$, where $K$ satisfies the H\"{o}rmander condition. 
Then the couple $\left(S_1,S_p\right)$ is $K$-closed in $\left(L^1\left(\ell^2\right),L^p\left(\ell^2\right)\right)$, with a constant depending only on the constant in H\"{o}rmander condition, the $L^2\left(\ell^2\right)\to L^2\left(\ell^2\right)$ norm of $P$ and the space $\left(X,d,\mu\right)$.\end{lem}
It has been originally stated in a scalar-valued version, but the proof translates verbatim into the $\ell^2$-valued version. We will use it only in the case of $p=2$ and $P$ given by the principal value.
\begin{cor}\label{kclt}The couple $\left(H^1(\T),H^2(\T)\right)$ is $K$-closed in $\left(L^1(\T),L^2(\T)\right)$.\end{cor}
\textit{Proof.} We apply Lemma \ref{bourgkcl} to the case of $X=\T$ with Lebesgue measure and metric $d(x,y)=\left|\e^{ix}-\e^{iy}\right|$, $S=\mathrm{span}\left\{\mathrm{e}^{int}:n\in\mathbb{N}\right\}$. The projection $P$ is just the Riesz projection corresponding to the kernel $K\left(x,y\right)=\frac{1}{1-\mathrm{e}^{i(x-y)}}$ which is known to satisfy the required conditions.\par 
From now on, our space of homogeneous type will be the Cantor group (with dyadic metric and Haar measure). The following trivial observation has to be made.
\begin{lem}\label{radial}If the convolution kernel $K$ is radial (i.e. $K(x)$ depends only on $d\left(x,\emptyset\right)$, or, equivalently, on $\min x$), then the integral in the H\"ormander condition vanishes.\end{lem}
\textit{Proof.} Suppose $d(x,\emptyset)\geq 2r$ and $d(y,\emptyset)<r$. Let $m$ be an integer such that $2^{-m-1}\leq r< 2^{-m}$. Then $d(x,\emptyset)\geq 2^{-m}$, $d(y,\emptyset)<2^{-m}$. Thus $\min x\leq m<\min y$, so $\min(x+y)=\min(x)$ as desired.\par 
Since $\left(\E_n-\E_{n-1}\right)f=\sum_{\max A=n}\left\langle f,w_A\right\rangle w_A$, convolution with the bounded kernel 
\[\kappa_n=\sum_{\max A=n}w_A\]
is exactly the operator $\Delta_n$. An easy calculation shows that 
\[\kappa_n(x)=\left\{\begin{array}{lll}1 &\ \mathrm{if}\ &n=0\\ 2^{n-1}&\ \mathrm{if}\ &1\leq n<\min x\\ -2^{n-1}&\ \mathrm{if}\ & n=\min x\\ 0&\ \mathrm{if}\ &n>\min x.\end{array}\right.\]
Sums over $n$ involving $\kappa_n$ are defined at each point $x\neq\emptyset$, since they automatically truncate at $n=\min x$. For $n\in\mathbb{N}$, we will denote by $e_n^*$ the $n$-th coordinate functional on $\ell^2$. \par 
A simple example of utilising Theorem \ref{czo} is an immediate proof of Theorem \ref{classicstein} in the dyadic case.
\begin{cor}The Stein martinagle inequality is true for the reversed dyadic filtration $\left(\mathcal{F}_{n-k}\right)_{k=0}^n$ on $[0,1]$ (or, equivalently, for the filtration $\left(\mathcal{F}^*_k\right)_{k=0}^n$).\end{cor}
\textit{Proof.} The equivalence of considering the inequality for the mentioned filtrations is achieved by changing the order of coordinates in the Cantor group, so we will prove it for the reversed dyadic filtration. The orthogonal projection $\left(f_k\right)_{k=0}^\infty\mapsto \left(\E_0 f_0,\ldots, \E_n f_n,0\ldots\right)$ is the convolution with the kernel \[\sum_{k\leq n}\left(\sum_{j\leq k}\kappa_j\right)\otimes \left(e_k\otimes e^*_k\right),\] since $\sum_{j\leq k}\kappa_j= \sum_{A\subset\{1,\ldots,k\}}w_A$ and $e_k\otimes e_k^*$ is the projection onto $k$-th coordinate. This kernel is expressed (finitely) in terms of kernels $\kappa_j$, so it is radial (and bounded), so Theorem \ref{czo} ends the proof.\par 
The following lemma will be used to transform limits arising in principal values into pointwise limits of martingales.
\begin{lem}\label{trunc} Let $v_1,\ldots,v_n,\ldots$ be vectors in the unit ball of a Banach space $B$. Denote by $K_m$ the $\mathcal{L}\left(\ell^2,B\right)$-valued kernel \[K_m=\sum_{n\leq m}\kappa_n\otimes\left(v_n\otimes e_n^*\right).\]
Then for any $f\in L^1\left(\ell^2\right)$, \[\lim_{m\to \infty}\int_{B(\emptyset,2^{-m})}K_m(y)f(x-y)\d \mu(y)=0\] for almost every $x$. \end{lem}
\textit{Proof.} Denote $\ell^2\ni u_m=\left(2^{-m},2^{-m+1},\ldots,1,0,\ldots\right)$ and $L^1\left(\ell^2\right)\ni |f|=\left(\left|f_n\right|\right)_{n=0}^\infty$. We have 
\[\aligned \left\|\int_{B(\emptyset,2^{-m})}K_m(y)f(x-y)\d \mu(y)\right\|_{B} &= \left\|\sum_{n\leq m}v_n\int_{B(\emptyset,2^{-m})}\kappa_n(y)f_n(x-y)\d \mu(y)\right\|_{B}\\ &\leq \sum_{n\leq m}2^{n-1}\int_{B(\emptyset,2^{-m})}\left|f_n(x-y)\right|\d \mu(y)\\ &= \sum_{n\leq m}\frac{2^{n-m-1}}{\mu\left(B(\emptyset,2^{-m})\right)}\int_{B(\emptyset,2^{-m})}\left|f_n(x-y)\right|\d \mu(y)\\ &=\sum_{n\leq m}2^{n-m-1}\E_{m}\left|f\right|(x)\\ &=\frac{1}{2}\left\langle u_m,\E_{m}|f|(x)\right\rangle\endaligned\]
which tends to $0$ for almost every $x$ since $\E_{m}|f|\to |f|$ almost everywhere and $u_m\to 0$ weakly in $\ell^2$. \par 
We are now ready to prove two main theorems of this section.
\begin{thm}\label{wktpl2}Let \[K=\sum_{n}\kappa_n\otimes\left(e_n\otimes e_n^*\right).\] The operator $f\mapsto \pv \left(K\ast f\right)$ coincides with the orthogonal projection $P:L^2\left(\ell^2\right)\to \iota\left(H^2(\delta)\right)$ for $f\in L^2\left(\ell^2\right)$ and is well defined and continuous as an $L^1\left(\ell^2\right)\to L^{1,\infty}\left(\ell^2\right)$ operator (and thus as an $L^p\left(\ell^2\right)\to L^p\left(\ell^2\right)$ one for $1<p<\infty$). \end{thm}
\textit{Proof.} Just as in the proof of Lemma \ref{trunc}, denote \[L^\infty\left(\mathcal{L}\left(\ell^2\right)\right)\ni K_m=\sum_{n\leq m} \kappa_n\otimes\left(e_n\otimes e_n^*\right).\] By the explicit form of $\kappa_n$, we have $K(y)=K_m(y)$ if $d(y,\emptyset)\geq 2^{-m}$. Let $f\in L^1\left(\ell^2\right)$. Then by Lemma \ref{trunc}
\[\aligned \pv \left(K\ast f\right)(x)&= \lim_{r\to 0}\int_{\mathbb{Z}_2^{\mathbb{N}^+}\setminus B\left(\emptyset,r\right)}K(y)f(x-y)\d \mu(y)\\
&= \lim_{m\to \infty}\int_{\mathbb{Z}_2^{\mathbb{N}^+}\setminus B\left(\emptyset,2^{-m}\right)}K(y)f(x-y)\d \mu(y)\\
&= \lim_{m\to \infty}\int_{\mathbb{Z}_2^{\mathbb{N}^+}\setminus B\left(\emptyset,2^{-m}\right)}K_m(y)f(x-y)\d \mu(y)\\ 
&= \lim_{m\to \infty}\int_{\mathbb{Z}_2^{\mathbb{N}^+}}K_m(y)f(x-y)\d \mu(y)\\
&=\lim_{m\to \infty}\left(K_m\ast f\right)(x).\endaligned\]
Let $f\in L^2\left(\ell^2\right)$. Since $e_n\otimes e_n^*$ is the projection onto $n$-th coordinate and $\kappa_n\ast f_n= \Delta_n f_n$, we have $K_m\ast f=\sum_{n\leq m}\left(\Delta_n f_n\right)\otimes e_n= \E_m\left(Pf\right)$, which tends to $Pf$, as desired. Now let $f\in L^1\left(\ell^2\right)$. The bounded kernel $K_m$ induces a norm 1 orthogonal projection and is radial, so by Lemma \ref{radial} and by Theorem \ref{czo}, \[\lambda\mu\left\{\left(\sum_{n\leq m}\left|\Delta_nf_n(x)\right|^2\right)^\frac{1}{2}>\lambda\right\}\leq C\left\|f\right\|_{L^1\left(\ell^2\right)}\]
where $C$ is a numerical constant. Passing to the limit with $m\to\infty$ we infer that $\sum_n \left|\Delta_nf_n(x)\right|^2$ is finite a.e., so the pointwise limit $\pv\left(K\ast f\right)(x)=\lim_{m\to\infty}\left(K_m\ast f\right)(x)=\lim_{m\to\infty}\left(\Delta_0f_0(x),\ldots,\Delta_m f_m(x),0,\ldots\right)$ exists almost everywhere and the result follows.\par
\begin{cor}\label{kcld}The couple $\left(\left(\iota\left(H^1(\delta)\right),\|\cdot\|_{L^1\left(\ell^2\right)}\right),\left(\iota\left(H^2(\delta)\right),\|\cdot\|_{L^2\left(\ell^2\right)}\right)\right)$ is $K$-closed in $\left(L^1\left(\ell^2\right), L^2\left(\ell^2\right)\right)$.\end{cor}
\textit{Proof.} It is simply a combination of Lemma \ref{bourgkcl} and Theorem \ref{wktpl2}. \par 

It has to be noted that $P$ is not $L^1\left(\ell^2\right)\to L^1\left(\ell^2\right)$ bounded. In fact, it is not even bounded when restriced to the subspace of $f$ such that $\E_{k-1}f_k=0$, which follows from Theorem \ref{drzewa} by restricting the domain to $\mathcal{F}_n$-measurable functions and reversing the order of the filtration. Consequently, the map $\iota^{-1}P$ can not be $L^1\left(\ell^2\right)\to L^1$ bounded, since that would imply \[\left\|\sum_n \pm \Delta_n f_n\right\|_{L^1}\lesssim \left\|f\right\|_{L^1\left(\ell^2\right)},\] which after averaging over choces of signs would give boundedness of $P$. However, we can prove a result analogous to Theorem \ref{wktpl2} for $\iota^{-1}P$. 

\begin{thm}\label{wktpsc}Let \[k=\sum_{n}\kappa_n\otimes e_n^*.\] The operator $f\mapsto \pv \left(k\ast f\right)$ coincides with the map $\iota^{-1}P:L^2\left(\ell^2\right)\to H^2(\delta)$ for $f\in L^2\left(\ell^2\right)$ and is well defined and continuous as an $L^1\left(\ell^2\right)\to L^{1,\infty}$ operator (and thus as an $L^p\left(\ell^2\right)\to L^p$ one for $1<p<\infty$).\end{thm}
\textit{Proof.} Analogously to the proof of Theorem \ref{wktpl2}, we put\[L^\infty\left(\left(\ell^2\right)^*\right)\ni k_m=\sum_{n\leq m}\kappa_n\otimes e_n^*,\] note that $\pv\left(k\ast f\right)(x)=\lim_{m\to\infty}k_m\ast f(x)$ for $f\in L^1\left(\ell^2\right)$ and that for $f\in L^2\left(\ell^2\right)$ we have $k_m\ast f=\E_m\left(\iota^{-1}Pf\right)\to \iota^{-1}Pf$. The only difference is the proof of existence of a pointwise limit of $k_m\ast f$ for $f\in L^1\left(\ell^2\right)$. \par 
By Theorem \ref{czo} applied to $k_m$, we get \[\left\|\sum_{a\leq n\leq b}\Delta_n f_n\right\|_{L^{1,\infty}}=\sup_{\lambda>0}\lambda\mu\left\{\left|\sum_{a\leq n\leq b}\Delta_n f_n\right|>\lambda\right\}\leq C\left\|\left(\sum_{a\leq n\leq b}\left|f_n\right|^2\right)^\frac{1}{2}\right\|_{L^1\left(\ell^2\right)}\]
which implies that $M_n=k_n\ast f=\sum_{j\leq n}\Delta_j f_j$ is a Cauchy sequence in $L^{1,\infty}$. Let $M$ be its limit in $L^{1,\infty}$ quasinorm. We are going to prove that $M_n$ converges to $M$ almost everywhere. By the weak type convergence, $M_n$ converges to $M$ in measure, so $M_{a_n}\to M$ almost everywhere for some increasing sequence $\left(a_n\right)$. Denote $I_n=\left\{a_n,\ldots,a_{n+1}-1\right\}$ and for any $j$ let $n(j)$ be the unique integer such that $j\in I_{n(j)}$. Consider the kernel 
\[L^{\infty}\left( \mathcal{L}\left(\ell^2,\left(\bigoplus_{n}\ell^\infty_{I_n}\right)_{\ell^2}\right)\right)\ni K_{\left(a_1,\ldots,a_m\right)}=\sum_{j<a_{m+1}}\kappa_j\otimes\left(v_j\otimes e_j^*\right),\]
where $\left(v_j\right)_{n}=0$ for $n\neq n(j)$ and \[\left(v_j\right)_{n(j)}(i)=\left\{\begin{array}{ll}1&\ \mathrm{if}\ i\geq j\ \mathrm{and}\ j\neq a_{n(j)}\\ 0&\ \mathrm{otherwise}\end{array}\right.\] for $i\in I_{n(j)}$. It is easy to see that $\left(K_{\left(a_1,\ldots,a_m\right)}\ast f\right)_{n}(j)=\sum_{a_{n}<i\leq j}\Delta_i f_i$ for $j\in I_n$. If $f\in L^2\left(\ell^2\right)$, then by the maximal inequality of Doob applied to the martingale $\left(\sum_{a_{n}<i\leq j}\Delta_i f_i\right)_{j\in I_n}$ with respect to the filtration $\left(\mathcal{F}_{j}\right)_{j\in I_n}$ we get \[\E \left(\max_{j\in I_n} \left|\sum_{a_{n}<i\leq j}\Delta_i f_i\right|\right)^2\lesssim \E\left|\sum_{a_n<i<a_{n+1}}\Delta_i f_i\right|^2\leq \E\sum_{a_n<i<a_{n+1}}\left|f_i\right|^2\]
which after summing over $n$ gives the $L^2$-boundedness of convolution with $K_{\left(a_1,\ldots,a_m\right)}$. Together with the radiality of $K_{\left(a_1,\ldots,a_m\right)}$, this gives the weak type estimate \[\sup_{\lambda}\lambda\mu\left\{\left(\sum_{n\leq m} \left(\max_{j\in I_n} \left|\sum_{a_{n}<i\leq j}\Delta_i f_i\right|\right)^2\right)^\frac{1}{2}>\lambda\right\}\lesssim \left\|f\right\|_{L^1\left(\ell^2\right)}.\] Thus $\sum_{n=1}^\infty \left(\max_{j\in I_n} \left|\sum_{a_{n}<i\leq j}\Delta_i f_i\right|\right)^2$ is finite a.e., so $\max_{j\in I_n} \left|\sum_{a_{n}<i\leq j}\Delta_i f_i\right|$ converges to $0$ a.e. But $\sum_{a_{n}<i\leq j}\Delta_i f_i=M_j-M_{a_n}$ for $j\in I_n$, so $M_j-M_{a_{n(j)}}$ converges to $0$, which together with $M_{a_{n(j)}}\to M$ completes the proof.

\section{Isomorphisms between independent sums}

We will use the identification of independent sums with certain subspaces of Orlicz spaces \cite{Dilworth}, \cite{JSch}.
\begin{deff} Let $\Phi:\left[0,\infty\right)\to \left[0,\infty\right)$ be convex, increasing and $\Phi(0)=0$. For $f\in L^0\left(X,B\right)$, where $X$ is a measure space and $B$ is a Banach space, the Orlicz norm associated with the function $\Phi$ is \[\|f\|_{\Phi}=\inf\left\{k>0: \int_X \Phi\left(\frac{\|f(t)\|_B}{k}\right)\d \mu\leq 1\right\}.\]\end{deff}
\begin{thm}\label{JSch} Let $\left(f_n\right)_{n=1}^\infty\in\left(\bigoplus_{n=1}^\infty L^1(\Omega,B)\right)_{\mathrm{ind}}$. Then \[\left\|\left(f_n\right)_{n=1}^\infty\right\|_{\mathrm{ind}}\sim \left\|\coprod_{n=1}^\infty f_{n}\right\|_{\Phi_1},\] where $\coprod_{n=1}^\infty f_{n}$ a function defined on the disjoint sum $\coprod_{n=1}^\infty \Omega_n $ such that $\left(\coprod_{n=1}^\infty f_n\right)\left(\omega\right)= f_k\left(\omega\right)$ whenever $\omega\in\Omega_k$ and \[\Phi_1(x)=\left\{\begin{array}{lll}x^2&\text{ for }&0\leq x\leq 1\\2x-1&\text{ for }&x>1.\end{array}\right.\]\end{thm}
Note the following routine but useful observations.
\begin{lem} \label{whatever}For $f\in L^0\left(X,B\right)$, \[\|f\|_{\Phi_1}\sim\inf_{\substack{g+h=f\\g\in L^1(B)\\ h\in L^2(B)}}\|g\|_{L^1(B)}+\|h\|_{L^2(B)}.\]\end{lem}
\textit{Proof.} To prove the inequality $\lesssim$, take any $g,h$ such that $g+h=f$. Then, by\linebreak $\Phi_1(x)\leq \min(2x,x^2)$ and convexity of $\Phi_1$,
\[\aligned \int \Phi_1\left(\frac{\|g+h\|_B}{\|g\|_{L^1(B)}+\|h\|_{L^2(B)}}\right)\d \mu&\leq \int \Phi_1\left(\frac{\|g\|_{L^1(B)}}{\|g\|_{L^1(B)}+\|h\|_{L^2(B)}}\cdot \frac{\|g\|_B}{\|g\|_{L^1(B)}}+\right.\\
&\left. \frac{\|h\|_{L^2(B)}}{\|g\|_{L^1(B)}+\|h\|_{L^2(B)}} \cdot\frac{\|h\|_B}{\|h\|_{L^2(B)}}\right)\d\mu\\ 
&\leq \frac{\|g\|_{L^1(B)}}{\|g\|_{L^1(B)}+\|h\|_{L^2(B)}}\int\Phi_1\left(\frac{\|g\|_B}{\|g\|_{L^1(B)}}\right)\d\mu\\
&+ \frac{\|h\|_{L^2(B)}}{\|g\|_{L^1(B)}+\|h\|_{L^2(B)}}\int\Phi_1\left(\frac{\|h\|_B}{\|h\|_{L^2(B)}}\right)\d\mu\\ &\leq 2.
 \endaligned\]
Since $\Phi_1\left(\frac{x}{2}\right)\leq \frac{\Phi_1(x)}{2}$, taking $k=2\left(\|g\|_{L^1(B)}+\|h\|_{L^2(B)}\right)$ gives the desired inequality. \par 
Now take $k$ such that $\int \Phi_1\left(\frac{\|f\|_B}{k}\right)\d\mu\leq 1$. We have to prove that $\inf_{g+h=f}\|g\|_{L^1(B)}+\|h\|_{L^2(B)}\lesssim k$. By homogenity we can assume that $k=1$. Let $g=f\cdot \mathbbm{1}_{\|f\|_B>1}$ and $h=f\cdot \mathbbm{1}_{\|f\|_B\leq 1}$. Then, by the inequality $x^\frac12\leq 1+x$, 
\[\aligned \|g\|_{L^1(B)}+\|h\|_{L^2(B)}&= \int_{\|f\|_B>1}\|f\|_B\d\mu+\left(\int_{\|f\|_B\leq 1}\|f\|_B^2\d\mu\right)^\frac12\\ &\leq 1+\int_{\|f\|_B>1}\left(2\|f\|_B-1\right)\d\mu+\int_{\|f\|_B\leq 1}\|f\|_B^2\d\mu\\&=1+\int\Phi_1\left(\|f\|_B\right)\d\mu\\&\lesssim 1.\endaligned\]
\begin{cor}\label{l1+l2}For $f\in\left(\bigoplus_{n=1}^\infty L^1(\Omega,B)\right)$ we have 
\[\indnorm{f}\sim \inf_{\substack{g_n+h_n=f_n\\ g_n\in L^1(\Omega,B)\\ h_n\in L^2(\Omega,B)}}\sum\|g_n\|_{L^1(\Omega,B)}+\sqrt{\sum \|h_n\|_{L^2(\Omega,B)}^2}.\]\end{cor}
\textit{Proof.} It is an immediate consequence of Theorem \ref{JSch} and Lemma \ref{whatever}.\par

Below we will be dealing with subspaces $X\subset L^1$ (where $L^1$ is possibly vector valued; we will sometimes omit indicating the space of values) such that $X\cap L^2$ is dense in $X$. Let $X,Y$ be two such subspaces, possibly with values in different Banach spaces. For an operator $T:X\to Y$, we denote $\|T\|_{\mathcal{L}\left(L^1\right)\cap \mathcal{L}\left(L^2\right)}=\|T\|_{\mathcal{L}\left(\left(X,\|\cdot\|_{L^1}\right), \left(Y,\|\cdot\|_{L^1}\right)\right)}+\|T\|_{\mathcal{L}\left(\left(X\cap L^2,\|\cdot\|_{L^2}\right), \left(Y\cap L^2,\|\cdot\|_{L^2}\right)\right)}$ (the sum of the operator norms induced by $L^1\to L^1$ and $L^2\to L^2$ norms).\par 

\begin{deff}Let $X_n\subset L^1$, $Y_n\subset L^1$ and consider a sequence of operators $\left(T_n\right)$, where $T_n:X_n\to Y_n$. We will call it $R_{\mathrm{ind}}$-bounded if the operator $\left(\bigoplus T_n\right)_{\mathrm{ind}}:\left(\bigoplus X_n\right)_{\mathrm{ind}}\to \left(\bigoplus Y_n\right)_{\mathrm{ind}}$ defined by \[\left(\bigoplus T_n\right)_{\mathrm{ind}}\left(f_n\right)=\left(T_n f_n\right)\] is continuous (or, in other words, $\indnorm{\left(T_n f_n\right)}\lesssim \indnorm{\left(f_n\right)}$).\end{deff}
The classical notion of $R$-boundedness arises when we consider operators defined on the whole $L^1$ space and replace $\indnorm{\cdot}$ with $\Lnorm{\cdot}$. \par

\begin{lem}\label{oplust} Let $T_n:X_n\to Y_n$ be operators between subspaces of $L^1$ spaces. Suppose that one of the following is satisfied.\\
(i) For any $n$ the couple $\left(\left(X_n,\|\cdot\|_{L^1}\right),\left(X_n\cap L^2,\|\cdot\|_{L^2}\right)\right)$ is $K$-closed in $\left(L^1,L^2\right)$ with a constant not dependent on $n$ and the operators are uniformly bounded in $L^1\to L^1$ and $L^2\to L^2$ norms. \\
(ii) The operators $T_n$ are uniformly bounded in the $L^1\to L^2$ norm.\\
Then the operators $T_n$ are $R_{\mathrm{ind}}$-bounded.
 \end{lem}
\textit{Proof.} Suppose (i) is satisfied. By Corollary \ref{l1+l2}
\[\aligned \left\|\left(T_n f_n\right)\right\|_{\mathrm{ind}}&\sim \inf_{u_n+v_n=T_n f_n}\sum\left\|u_n\right\|_{L^1}+ \sqrt{\sum\left\|v_n\right\|^2_{L^2}}\\
&\leq \inf_{\substack{g_n+h_n=f_n\\ g_n\in X_n,h_n\in X_n\cap L^2}}\sum\left\|T_n g_n\right\|_{L^1}+ \left(\sum \left\|T_n h_n\right\|_{L^2}^2\right)^\frac{1}{2}\\ 
&\lesssim \inf_{\substack{g_n+h_n=f_n\\ g_n\in X_n,h_n\in X_n\cap L^2}}\sum\left\|g_n\right\|_{L^1}+ \left(\sum \left\|h_n\right\|_{L^2}^2\right)^\frac{1}{2}\\
&\lesssim \inf_{\substack{g_n+h_n=f_n\\ g_n\in L^1,h_n\in L^2}}\sum\left\|g_n\right\|_{L^1}+ \left(\sum \left\|h_n\right\|_{L^2}^2\right)^\frac{1}{2}\\
&\sim\left\|\left(f_n\right)\right\|_{\mathrm{ind}},\endaligned\] 
where the first '$\lesssim$' inequality follows from the uniform bound on norms of $T_n$ and the second from the uniform bound for the $K$-closedness constants and Lemma \ref{kcleq}.   \par 
Suppose now that (ii) is satisfied. Take any decomposition $f_n=g_n+h_n$. Then 
\[\indnorm{\left(T_nf_n\right)}\lesssim \sqrt{\sum\left\|T_n f_n\right\|_{L^2}^2}\lesssim \sqrt{\sum\left\|f_n\right\|_{L^1}^2}\]\[\leq \sqrt{\sum\left\|u_n\right\|_{L^1}^2}+ \sqrt{\sum\left\|v_n\right\|_{L^1}^2}\leq \sum\left\|u_n\right\|_{L^1}+ \sqrt{\sum\left\|v_n\right\|_{L^2}^2}\] 
and taking the infimum over all choices of $g_n,h_n$ completes the proof.\par 
It has to be noted that the $K$-closedness assumption in (i) can not be omitted, since one can find, which we will not do, sequences of functions such that $\left\|F_n\right\|_p=\left\|f_n\right\|_p$ for $p\in\{1,2\}$ and $\indnorm{\left(F_n\right)} \indnorm{\left(f_n\right)}^{-1}$ is arbitrarily large. Nonetheless, the $K$-closedness assumption is trivially satified if $X_n=L^1$. \par 
We will return for a moment to the issues of Section 2. Namely, we prove a quantitative version of Corollary \ref{renorm}. 
\begin{lem}\label{schodkapr} Let $1\leq p\leq\infty$, $N\in\mathbb{N}^+$ and for $f\in W^{1,p}[0,1]$ denote $(\E_N f)(x):=N\int_{I_N(x)}f(y)\d y$, where $I_N(x)$ is an interval of the form $\left[\frac{j}{N},\frac{j+1}{N}\right)$ such that $x\in I_N(x)$. Then \[\left\|\left(\mathrm{id}-\E_N\right)f\right\|_{L^p}\leq \frac{1}{N}\left\|f'\right\|_{L^p}.\]\end{lem}
\textit{Proof.} For any $x$ we have \[\left|\left(\mathrm{id}-\E_N\right)f(x)\right|=\left|N\int_{I_N(x)}\left(f(x)-f(y)\right)\d y\right|= \left|N\int_{I_N(x)}\int_{y}^x f'(z)\d z\d y\right|\leq\]\[\leq N\int_{\left(I_N(x)\right)^2}\left|f'(z)\right| \mathbbm{1}_{\mathrm{conv}\{x,y\}}(z)\d z \d y\leq \int_{I_N(x)}\left|f'(z)\right|\d z.\]
Thus \[\int_0^1 \left|\left(\mathrm{id}-\E_N\right)f(x)\right|^p\d x= \sum_j \int_{\frac{j}{N}}^\frac{j+1}{N}\left|\left(\mathrm{id}-\E_N\right)f(x)\right|^p\d x\leq \]\[\leq \sum_j \frac{1}{N}\left(\int_{\frac{j}{N}}^\frac{j+1}{N}\left|f'(z)\right|\d z\right)^p\leq \frac{1}{N}\sum_j \frac{1}{N^{p-1}} \int_{\frac{j}{N}}^\frac{j+1}{N}\left|f'(z)\right|^p\d z= \frac{1}{N^p}\int_0^1 \left|f'(z)\right|^p\]
as desired.
\begin{thm}\label{transference}Let $N_1,N_2,\ldots$ be positive integers such that $N_k\mid N_{k+1}$ for any $k$. Denote by $\mathcal{F}_k$ the sigma-algebra generated by intervals of the form $\left[\frac{2\pi j}{N_k},\frac{2\pi(j+1)}{N_k}\right)$ for $0\leq j<N_k$ and by $\mathcal{F}^*_k$ the sigma-algebra of subsets of $\T$ whose chracteristic functions are $\frac{2\pi}{N_k}$-periodic. Suppose that $f_1,f_2,\ldots$ are trigonometric polynomials such that $\mathrm{deg}f_k\leq d_k$ and $f_k$ is $\mathcal{F}^*_k$-measurable. Then \[\left(1-C\theta\right)\left\|\left(f_k\right)\right\|_{\mathrm{ind}}\leq \left\|\left(f_k\right)\right\|_{L^1\left(\ell^2\right)}\leq \left(1+C\theta\right)\left\|\left(f_k\right)\right\|_{\mathrm{ind}},\]
where $\theta=\sup_k \frac{d_k}{N_{k+1}}$ and $C$ is a numerical constant. \end{thm}
\textit{Proof.} Note that $\mathcal{F}_k$ and $\mathcal{F}_k^*$ are independent and for general reasons, $\E_k^* \E_{k+1}=\E_{k+1}\E^*_k$. For any $k$, the function $\E_{k+1}f_k$ is independent of the sigma-algebra generated by the functions $\E_{j+1}f_j$, where $j<k$. Indeed, they are $\mathcal{F}_{j+1}$-measurable and thus $\mathcal{F}_k$-measurable, while $\E_{k+1}f_k=\E_{k+1}\E^*_{k}f_k=\E^*_k\E_{k+1}f_k$ is $\mathcal{F}_k^*$-measurable.  Thus, the functions $\E_{k+1}f_k$ are independent random variables. By Corrolary \ref{2indgeqL1l2},
\[\left|\left\|\left(f_k\right)\right\|_{L^1\left(\ell^2\right)}- \left\|\left(f_k\right)\right\|_{\mathrm{ind}}\right|\leq \left|\left\|\left(\E_{k+1}f_k\right)\right\|_{L^1\left(\ell^2\right)}- \left\|\left(\E_{k+1}f_k\right)\right\|_{\mathrm{ind}}\right|+ \]\[+\left| \left\|\left(f_k\right)\right\|_{L^1\left(\ell^2\right)} - \left\|\left(\E_{k+1}f_k\right)\right\|_{L^1\left(\ell^2\right)}\right|+ \left|\left\|\left(f_k\right)\right\|_{\mathrm{ind}}- \left\|\left(\E_{k+1}f_k\right)\right\|_{\mathrm{ind}}\right|\lesssim \]\[\lesssim \left\|\left(\left(\mathrm{id}-\E_{k+1}\right)f_k\right)\right\|_{\mathrm{ind}} =\left\|\left(\left(\mathrm{id}-\E_{k+1}\right)V_{d_k}f_k\right)\right\|_{\mathrm{ind}},\]
where $V_{d_k}$ is the de la Vall\'{e}e-Poussin kernel of order $d_k$. By Lemma \ref{schodkapr}, for $p\in\{1,2\}$ we have $\left\|\left(\mathrm{id}-\E_{k+1}\right)V_{d_k}f\right\|_{L^p}\leq \frac{2\pi }{N_{k+1}}\left\|\left(V_{d_k}f\right)'\right\|_{L^p}\leq \frac{2\pi d_k}{N_{k+1}}\left\|V_{d_k}f\right\|_{L^q}\lesssim\frac{d_k}{N_{k+1}}\left\|f\right\|_{L^q}$ due to Bernstein inequality. Application of Lemma \ref{oplust} to operators $\left(\mathrm{id}-\E_{k+1}\right)V_{d_k}$ acting on the whole $L^1(\T)$ ends the proof.

\begin{cor}\label{crudetransference}Let $\left(f_1,\ldots,f_n,\ldots\right)\in L^1\left(\T,\ell^2(B)\right)$ and let $T$ be dilation $Tf(x)=f\left(2x (\mathrm{mod}\ 2\pi)\right)$. Then \[\lim_{a\to \infty}\left\|\left(T^{ak}f_k\right)\right\|_{L^1\left(\ell^2(B)\right)}= \indnorm{\left(f_k\right)}.\] Trivially, the same holds if we replace $\T$ with $[0,1]$ and $2\pi$ with $1$. \end{cor}
\textit{Proof.} It is enough to handle the scalar-valued case since we can consider $\left\|f_k\right\|_B$ instead of $f_k$. Let $\varepsilon>0$. Take polynomials $g_1,\ldots,g_k$ such that $\indnorm{\left(f_k-g_k\right)}<\varepsilon$. Then $T^{ak}g_k$ is of degree $2^{ak}\deg g_k$ and is $\frac{2\pi}{2^{ak}}$-periodic. For some global constant $C$, by Corollary \ref{2indgeqL1l2} and Theorem \ref{transference}
\[\left| \left\|\left(T^{ak}f_k\right)\right\|_{L^1\left(\ell^2\right)}- \indnorm{\left(f_k\right)}\right|\leq \left| \left\|\left(T^{ak}g_k\right)\right\|_{L^1\left(\ell^2\right)}- \indnorm{\left(g_k\right)}\right|+ \]\[\left|\indnorm{\left(f_k\right)}- \indnorm{\left(g_k\right)}\right|+ \left|\left\|\left(T^{ak}f_k\right)\right\|_{L^1\left(\ell^2\right)}- \left\|\left(T^{ak}g_k\right)\right\|_{L^1\left(\ell^2\right)}\right|\lesssim\] \[C \max_{k} \frac{2^{ak}}{2^{a(k+1)}}+C\varepsilon= C\left(2^{-a}+\varepsilon\right).\]
Thus, by a suitable choice of $a$ and $\varepsilon$ we can make $\left|\left\|\left(T^{ak}f_k\right)\right\|_{L^1\left(\ell^2\right)}-\indnorm{\left(f_k\right)}\right|$ as small as desired.

\begin{thm}\label{pelcz}Let $X,Y$ be subspaces of $L^1$ spaces, possibly on different domains and possibly taking values in different Banach spaces, such that there is an isomorphism $M:X\to Y$ such that the sequences $\left(M,M,\ldots \right)$ and $\left(M^{-1},M^{-1},\ldots\right)$ are $R_{\mathrm{ind}}$-bounded. We will say that a sequence $\left(X_n\right)_{n=1}^\infty$ of subspaces of $X$ is a good net if it satisfies the following conditions:\\
(i) $X_n\subset L^2$ are finite dimensional\\
(ii) there exist operators $P_n:X\to X_n$ and $i_n:X_n\to X$ such that $P_n i_n=\mathrm{id}_{X_n}$ and sequences $\left(P_n\right)$ and $\left(i_n\right)$ are indenpendently $R$-bounded \\
(iii) $\bigcup X_n$ is dense in $X$ in the $L^2$ norm\\
(iv) for any $n_1,n_2$ there is $n_3$ such that $X_{n_1}\cup X_{n_2}\subset X_{n_3}$.\\
Suppose that $\left(X_n\right)_{n=1}^\infty$ is a good net in $X$ and $\left(Y_n\right)_{n=1}^\infty$ is a good net in $Y$. Then \[\left(\bigoplus_{n=1}^\infty X_n\right)_{\mathrm{ind}}\sim \left(\bigoplus_{n=1}^\infty Y_n\right)_{\mathrm{ind}}.\]\end{thm}
\textit{Proof.} The core of the proof is the fact that $\left(\bigoplus X_n\right)_{\mathrm{ind}}$ is isomorphic to a complemented subspace of $\left(\bigoplus Y_n\right)_{\mathrm{ind}}$. Consider the finite dimensional subspace $Mi_n X_n\subset Y$. The projection from $Y$ onto this subsapce is given by the formula $Q_n=Mi_nP_nM^{-1}$. By the conditions (iii) and (iv) for $\left(Y_n\right)$, we can find $\tilde{n}$ such that $\left(\mathrm{id}+T_n\right)Mi_n X_n\subset Y_{\tilde{n}}$ and $\left\|T_n\right\|_{L^1\to L^2}$ is as small as desired.  Obviously $\left(\mathrm{id}+T_n\right)Mi_n X_n= \left(\mathrm{id}+T_nQ_n\right)Mi_n X_n$. By the (ii) part of Lemma \ref{oplust}, the sequences $\left(\mathrm{id}+T_nQ_n\right)$ and $\left(\left(\mathrm{id}+T_nQ_n\right)^{-1}\right)$ are $R_{\mathrm{ind}}$-bounded. Then $\left(\bigoplus X_n\right)_{\mathrm{ind}}$ is isomorphic to $\left(\bigoplus \left(\mathrm{id}+T_nQ_n\right)Mi_n X_n\right)_{\mathrm{ind}}\subset \left(\bigoplus Y_{\tilde{n}}\right)_{\mathrm{ind}}$ via the map $\left(\bigoplus \left(\mathrm{id}+T_nQ_n\right)Mi_n \right)_{\mathrm{ind}}$ (since the inverses of restrictions of $i_n$ to $i_n X_n$ coincide with $P_n$, which are $R_{\mathrm{ind}}$-bounded). The projection from $\left(\bigoplus Y_{\tilde{n}}\right)_{\mathrm{ind}}$, which is complemented in $\left(\bigoplus Y_{n}\right)_{\mathrm{ind}}$ onto $\left(\bigoplus \left(\mathrm{id}+T_nQ_n\right)Mi_n X_n\right)_{\mathrm{ind}}$ is $$\left(\bigoplus \left(\mathrm{id}+T_nQ_n\right)Q_n\left(\mathrm{id}+T_nQ_n\right)^{-1} \right)_{\mathrm{ind}}$$.

\par 
We finish the proof in the spirit of Pełczy\'{n}ski decomposition principle. Let $\sigma:\mathbb{N}\to\mathbb{N}$ be any function satisfying $\# \sigma^{-1}(n)=\aleph_0$ for any $n$ and denote $V=\left(\bigoplus X_{n}\right)_{\mathrm{ind}}$, $W=\left(\bigoplus Y_{n}\right)_{\mathrm{ind}}$, $\tilde{V}=\left(\bigoplus X_{\sigma(n)}\right)_{\mathrm{ind}}$, $\tilde{W}=\left(\bigoplus Y_{\sigma(n)}\right)_{\mathrm{ind}}$. Obviously $\left(X_{\sigma(n)}\right)$ is a good net in $X$ and $\left(Y_{\sigma(n)}\right)$ is a good net in $Y$. By the already proved part of the theorem, $\tilde{V}$ is complemented in $V$ and by obvious properties of the independent sum, $V\sim \tilde{V}\oplus A\sim \tilde{V}\oplus\tilde{V}\oplus A\sim \tilde{V}\oplus V\sim \tilde{V}$ for some $A$. By symmetry, $W\sim \tilde{W}$. But again, by the already proved part, $\tilde{V}=\tilde{W}\oplus B$ for some $B$. Thus $\tilde{V}\sim \tilde{W}\oplus B\sim \tilde{W}\oplus \tilde{W}\oplus B\sim \tilde{W}\oplus \tilde{V}$. By symmetry, $\tilde{W}\sim \tilde{V}\oplus\tilde{W}$ and ultimately $V\sim \tilde{V}\sim \tilde{V}\oplus \tilde{W}\sim\tilde{W}\sim W$. 
\par 
We recall a claasical Marcinkiewicz-Zygmund inequality. 
\begin{lem}\label{marzyg} Let $B_1,B_2$ be Hilbert spaces and $T:X\to L^1\left(B_2\right)$, where $X\subset L^1\left(B_1\right)$ be bounded. Then for any $f_1,f_2,\ldots\in L^1\left(B_1\right)$,
\[\int\left(\sum \left\|Tf_n(x)\right\|_{B_2}^2\right)^\frac{1}{2}\d x\leq \|T\| \int\left(\sum \left\|f_n(x)\right\|_{B_1}^2\right)^\frac{1}{2}\d x.\]\end{lem}

\begin{thm}Let $\limsup a_n=\limsup b_n=\infty$. Then the spaces $\left(\bigoplus  H^1_{a_n}(\T)\right)_{\mathrm{ind}}$ and $\left(\bigoplus \iota\left( H^1_{b_n}(\delta)\right)\right)_{\mathrm{ind}}$ are isomorphic. \end{thm}

\textit{Proof.} We want to apply Theorem \ref{pelcz} to the case $X=\iota\left(H^1(\delta)\right)$, $X_n=\iota\left(H^1_{b_n}(\delta)\right)$, $Y=H^1(\T)$, $Y_n=H^1_{a_n}(\T)$. There are two nontrivial facts we need to verify. First is the existence of isomorphism $M$ and the second is the condition (ii) in the trigonometric case (since in the dyadic case $i_n$ are just identities and $P_n=\E_n$ act on the whole $L^1\left(\ell^2\right)$ with uniformly bounded $L^1\to L^1$ and $L^2\to L^2$ norms).\par 
\textit{Approach 1.} We will make use of the $K$-closedness. The existence of the desired isomorphism $M$ is guaranteed by the following theorem of Wojtaszczyk \cite{smallwojt}.\par 
\begin{thm}\label{franklin} The $L^2$-normalised Haar system in $H^1(\delta)$ is equivalent to the orthonormal Franklin system in $H^1(\T)$. \end{thm}
Since $\iota$ is an isometry in $L^2$ norms, the isomorphism $M$ transforming the Franklin system into image of the Haar system in $\iota$ is an isometry in $L^2$ norms. Thus, by Lemma \ref{oplust} and Corollaries \ref{kclt}, \ref{kcld}, sequences $\left(M\right)_{n=1}^\infty$ and $\left(M^{-1}\right)_{n=1}^\infty$ are $R_{\mathrm{ind}}$-bounded. In order to prove condition (ii) for $H^1_n(\T)$, we will use an observation of Bourgain and Pełczyński \cite{bigwojt}. 
\begin{lem}\label{bourgpelcz}There exist uniformly isomorphic and uniformly complemented copies of $H^1_n(\T)$ in $H^1(\T)$, in the sense of $\mathcal{L}\left(L^1\right)\cap \mathcal{L}\left(L^2\right)$ norm.\end{lem}
\textit{Proof.} Denote by $K_n$ the Fej\'{e}r kernel and $f\oplus g=\sum \widehat{f}(n)\e^{i2nt}+\sum \widehat{g}(n)\e^{i(2n+1)t}$. Let $R_n:H^1_n(\T)\to H^1_n(\T)$, $i_n: H^1_n(\T)\to H^1(\T)$, $P_n:H^1(\T)\to H^1_n(\T)$ be defined by $R_n\left(\sum_{k=0}^n a_k \e^{\i k t}\right)=\sum_{k=0}^n a_{n-k}\e^{\i k t}$, $i_n(f)=f\oplus R_n f$, $P_n\left(f\oplus g\right)= K_n\ast f + R_n\left(K_n\ast g\right)$. One easily checks that $P_n i_n=\mathrm{id}_{H^1_n(\T)}$ and $P_n$, $i_n$ are uniformly bounded in any $L^p$ norm.\par
Since by Lemma \ref{kcleq} and Corollary \ref{kclt} the uniform complementation of $i_n\left(H^1_n(\T)\right)$ in $H^1(\T)$ gives $K$-closedness of $$\left(\left(i_n\left(H^1_n(\T)\right),\|\cdot\|_{L^1}\right), \left(i_n\left(H^1_n(\T)\right),\|\cdot\|_{L^2}\right)\right)$$ in $\left(L^1,L^2\right)$, the sequences $\left(i_n\right)$ and $\left(P_n\right)$ are $R_{\mathrm{ind}}$-bounded by Lemma \ref{oplust}.\par 
\textit{Approach 2.} We will not make use of the $K$-closedness results. The fact that sequences $\left(P_n\right)$ and $\left(i_n\right)$ from Lemma \ref{bourgpelcz} are $R_{\mathrm{ind}}$-bounded can be verified directly, by Lemma \ref{oplust} for operators acting on the whole $L^1$. Indeed, it is easy to see that $\left|R_nf\right|$ and $\left|f\right|$ have the same dirstribution and $\indnorm{\left(f_n\oplus g_n\right)}\sim \indnorm{\left(f_n\right)}+ \indnorm{\left(g_n\right)}$, so the result follows from uniform boundedness of the Fej\'{e}r kernels. In order to prove the $R_{\mathrm{ind}}$-boundedness of the isomorphism, we need a stronger version of Theorem \ref{franklin}, proved by Meyer in \cite{meyer}.
\begin{thm}\label{meywavelet}The $L^2$ normalised Haar system in $H^1(\delta)$ is equivalent to the Meyer wavelet system.\end{thm}
Denote by $M$ the isomorphism of $H^1(\T)\to \iota\left(H^1(\delta)\right)$ given by Theorem \ref{meywavelet} and by $T_{\T}$ and $T_{[0,1]}$ the dilations on $\T$ and $[0,1]$ respectively. We have $MT_{\T}=T_{[0,1]}M$, since it is enough to check this identity for the Haar basis, and it reduces to the fact that the wavelet system satisfies the same dilation equation as the Haar system does: $T_{[0,1]}h_{2^n+j}=2^{-\frac12}\left(h_{2^{n+1}+j}+ h_{2^{n+1}+2^n+j}\right)$. Therefore, by Lemma \ref{marzyg}, \[\left\|\left(T^{an}f_n\right)\right\|_{L^1\left(\ell^2\right)}\sim \left\|\left(MT^{an}f_n\right)\right\|_{L^1\left(\ell^2\left(\ell^2\right)\right)}= \left\|\left(T^{an}Mf_n\right)\right\|_{L^1\left(\ell^2\left(\ell^2\right)\right)}\] for any $a$, which by Corollary \ref{crudetransference} gives $\indnorm{\left(f_n\right)}\sim \indnorm{\left(Mf_n\right)}$ as desired.

\section{Remarks and open questions}
We note that if our only interest in Section 4 was the $K$-closedness result, any potential problems with existence of the principal value could have been remedied by considering bounded kernels corresponding to projections onto $H^1_n(\delta)$. However, our goal was to develop a principal value result parallel to the Riesz projection $L^1(\T)\to H^1(\T)$. \par 
One can define $\mathrm{Ber}_T\ f=\frac{\E \left|Tf\right|}{\E \left|f\right|}$ for a (non necessarily bounded) operator $T:X\rightarrow L^1\left(\mathbb{T}\right)$ and ask when Lemma \ref{ber} is true for $\|\cdot\|_{\ell^2_2}$. In $\mathbb{R}$, it is enough for $T$ to commute with integration and satisfy $\mathrm{supp}Tf\subset \mathrm{supp}f$. In $\mathbb{C}$ these conditions are not sufficient, as for example $Tf=f'-\i f$ fails the lemma.\par 
The most known version of the Pełczyński decomposition principle mentioned in the proof of Theorem \ref{pelcz} is as follows \cite{bigwojt}.
\begin{thm} Let $X,Y$ be two Banach spaces such that $X$ is complemented in $Y$ and $Y$ is complemented in $X$. Assume additionally that $X$ is infinitely divisible, i.e. $X\sim \left(\bigoplus_{n=1}^\infty X\right)_p$ for some $p$. Then $X\sim Y$. \end{thm}
An obvious modification of the proof allows to replace the assumption $X\sim \left(\bigoplus_{n=1}^\infty X\right)_p$ with $X\sim \left(\bigoplus_{n=1}^\infty X\right)_{\mathrm{ind}}$. This, together with the appearance of the space $\left(\bigoplus X_{\sigma(n)}\right)_{\mathrm{ind}}$ in the proof of Theorem \ref{pelcz}, would suggest that this space is infinitely divisible in the sense of independent sum. However, unlike in the case of direct $\ell^p$-sums, there is no apparent reason for independent sum to be associative in the sense \[\left(\bigoplus_{n}\left(\bigoplus_{i\in I_n} X_i\right)_{\mathrm{ind}}\right)_{\mathrm{ind}}\sim \left(\bigoplus_{i\in \bigcup I_n} X_i\right)_{\mathrm{ind}}\]
for a disjoint family $\left\{I_n:n\in\mathbb{N}\right\}$. 
\par
It has been proved by Maurey \cite{Maurey} that $H^1(\delta)$ and $H^1(\T)$ are isomorphic, but the proof did not give an explicit way to construct an unconditional basis in $H^1(\T)$. Carleson \cite{Carleson} gave an example of a system in $H^1(\T)$ equivalent to the Haar basis in $H^1(\delta)$. Shortly after, Wojtaszczyk modified his proof to achieve orthonormality \cite{smallwojt}. Not much is known about the isomorphism class of $\left(\bigoplus H^1_n(\T)\right)_{\mathrm{ind}}$. The Rosenthal space $R$ is the sequence space normed by \[\left\|\left(a_n\right)\right\|_{R}=\|\left(a_n\right)\|_{\infty}+ \|\left(a_n \sqrt{w_n}\right)\|_2\] where $\sum w_n=\infty$ and $w_n\to 0$ (different choices of $\left(w_n\right)$ give isomorphic spaces). Its predual is the sequence space with a norm equivalent to \[\left\|\left(a_n\right)\right\|_{R_*}= \left(\sum \min\left(\left|a_n\right|,a_n^2w_n^{-1}\right)\right)^\frac12.\] It is known that $R_*$ is complemented in $H^1(\delta)$ and $R$ is complemented in $\mathrm{BMO}(\delta)$. The difficulty in handling the space $\left(\bigoplus H^1_n(\T)\right)_{\mathrm{ind}}$ can be seen in the fact that is not known whether it is isomorphic to $R_*$. A more detailed discussion can be found in \cite{smallmuller} and \cite{bigmuller}. A curious thing is that despite the Maurey isomorphism, $\left(\bigoplus H^1_n(\T)\right)_{\mathrm{ind}}$ is isomorphic to a complemented and invariant subspace of in $H^1(\delta)$, while the natural (based on a minor modification of the proof of Theorem \ref{mnoznik}) invariant embedding of this space into $H^1(\T)$ is not complemented by the Klemes theorem. This shows a substantial difference between the structure of invariant and complemented subspaces of $H^1(\T)$ and $H^1(\T\times\T)$. The non-isomorphism of these spaces is a nontrivial fact \cite{Bourgain}. \par
By an argument identical to the proof of Lemma \ref{rudin}, one can prove that any projection $L^1\left(\ell^2\right)\to \iota\left(H^1(\delta)\right)$ has to be of the form \[f=\sum_{A}\widehat{f}(A)w_A\mapsto \sum_{A} \left\langle \widehat{f}(A),u_A\right\rangle e_{\max A} w_A\] with $\left\langle u_A, e_{\max A}\right\rangle=1$. However, it is not clear to us how to deduce that $u_A=\e_{\max A}$, giving an orthogonal projection unbounded by the remark preceding the proof of Theorem \ref{wktpsc}.

\end{document}